%% file: main.tex
\documentclass[a4paper,11pt,oneside]{article}

\input{header.tex}
\usepackage{siunitx}
\usepackage{lineno}
\usepackage{hyperref}
\usepackage{placeins}
\newcommand\plot[1]{\let\frame\relax
  \frame{\includegraphics[clip,trim=0 220 0 60,width=8cm]{#1}}}

\newcommand\tsb[1]{_{\textup{#1}}}
\newcommand\Ieff{I\tsb{eff}}

\usepackage{authblk}
\usepackage[top=2cm, bottom=2cm, left=2cm, right=2cm]{geometry}

\theoremstyle{plain} 
\newtheorem{theorem}{Theorem}[section]

\newtheorem{Lemma}[theorem]{Lemma}

\newtheorem{Remark}[theorem]{Remark}
\newtheorem{Definition}[theorem]{Definition}

\newtheorem{Problem}[theorem]{Problem}

\theoremstyle{definition} %

\newtheorem{remark}[theorem]{Remark}

\theoremstyle{remark} %

\usepackage{newunicodechar}
\newunicodechar{ℝ}{\mathbb{R}}
\newunicodechar{Ψ}{\Psi}
\newunicodechar{α}{\alpha}
\newunicodechar{∇}{\nabla}
\newunicodechar{σ}{\sigma}
\newunicodechar{ν}{\nu}
\newunicodechar{ρ}{\rho}
\newunicodechar{ϑ}{\vartheta}
\newunicodechar{Γ}{\Gamma}

\usepackage{xcolor}

\newcommand{\Q}[1]{Q_{#1}}

\makeatletter
\newcommand{\customlabel}[2]{%
  \protected@write \@auxout {}{\string \newlabel {#1}{{#2}{\thepage}{#2}{#1}{}} }%
  \hypertarget{#1}{#2}%
}
\makeatother
\usepackage{tikz}
\usetikzlibrary{matrix,calc}
\usepackage{pgfplots}
\usepgfplotslibrary{groupplots,colorbrewer}
\usepgfplotslibrary{colormaps} \pgfplotsset{ colormap/Set1-5, cycle multiindex*
  list={ mark list*\nextlist Set1-5\nextlist }, every axis/.append style =
  {thick},%
}%
\pgfplotsset{tick style = {thick,black}}%
\newcommand{\doi}[1]{\href{https://doi.org/#1}{doi:#1}}

\begin{document}

    \title{
      {Anisotropic space-time goal-oriented error control and mesh adaptivity
        for convection-diffusion-reaction equations}
    }

    \author[1]{M. Bause}
    \author[1]{M. P. Bruchh\"auser}
    \author[2,3]{B. Endtmayer}
    \author[1]{N. Margenberg}
    \author[4]{I. Toulopoulos}
    \author[2,3]{T. Wick}

    \affil[1]{Helmut Schmidt University, University of the German Federal Armed Forces Hamburg,
      Faculty of Mechanical and Civil Engineering, Chair of Numerical Mathematics,
      Holstenhofweg 85,
      22043 Hamburg,
      Germany}

    \affil[2]{Leibniz Universit\"at Hannover, Institut f\"ur Angewandte
      Mathematik, Welfengarten 1, 30167 Hannover, Germany}

    \affil[3]{Cluster of Excellence PhoenixD (Photonics, Optics, and
      Engineering -- Innovation Across Disciplines), Leibniz Universit\"at Hannover, Germany}

    \affil[4]{University of Western Macedonia, Department of Informatics, Fourka Area, 52100 Kastoria, Greece}

    \maketitle
\begin{abstract}
We present an anisotropic goal-oriented error estimator based on the Dual
Weighted Residual (DWR) method for time-dependent convection-diffusion-reaction (CDR) equations.
Using anisotropic interpolation operators the estimator is elementwise separated with respect to the single directions in space and time leading to adaptive, anisotropic mesh refinement in a natural way.
To prevent spurious oscillations the streamline upwind Petrov-Galerkin (SUPG) method
is applied to stabilize the underlying system in the case of high P\'{e}clet numbers.
Efficiency and robustness of the underlying algorithm are demonstrated for
different goal functionals.
The directional error indicators quantify anisotropy of the solution with
respect to the goal, and produce meshes that efficiently capture sharp layers.
Numerical examples show the superiority of the proposed approach over isotropic
adaptive and global mesh refinement using established benchmarks for
convection-dominated transport.
\end{abstract}


    \section{Introduction}
    Adaptive finite element methods (AFEM) are essential for solving partial
    differential equations with localized phenomena such as layers, shocks, or
    singularities. For applications involving such phenomena,
    convection-diffusion-reaction (CDR) equations have been established as a prototypical model for more
    sophisticated systems of interest in practice, for instance, the
    Navier-Stokes equations.
    Their numerical solution remains an active topic of research. Following John et
    al.~\cite{BruchhaeuserJKN18}, we regard a numerical solution \emph{adequate} if it \begin{enumerate*}[label=(\roman*)]
    \item[\customlabel{itm:cdr:1}{(O1)}] captures \emph{sharp layers},
    \item[\customlabel{itm:cdr:2}{(O2)}] does not exhibit \emph{spurious oscillations} and
    \item[\customlabel{itm:cdr:3}{(O3)}] is \emph{feasible} and \emph{efficient} to compute.
    \end{enumerate*}
    In the past, a lot of effort has been spent on addressing the
    objectives~\ref{itm:cdr:1} and~\ref{itm:cdr:2}. Typically, these issues are
    addressed by means of stabilization. For many years, the streamline upwind
    Petrov–Galerkin (SUPG) method~\cite{BruchhaeuserHB79,BruchhaeuserBH82} and
    other residual-based stabilizations have been the most widely used
    techniques. Another approach that has gained traction is the class of
    algebraically stabilized schemes, which relies on the algebraic system from
    the Galerkin finite element discretization and introduce 
    limiters for controlling the slopes of the concrete solution.
    For a review of flux corrected transport we
    refer to~\cite{kuzmin2024property,kuzmin_flux-corrected_2012} and references
    therein. We refer
    to~\cite{BruchhaeuserRST08,augustin_assessment_2011,tobiska_robust_2015,BruchhaeuserJKN18} for a general review of stabilization techniques. For a numerical comparison
    related to time-dependent CDR equations we refer
    to~\cite{BruchhaeuserJS08,codina_comparison_1998}. In this work we
    address~\ref{itm:cdr:1} and~\ref{itm:cdr:2} through SUPG stabilization.
    We aim to address~\ref{itm:cdr:3} by means of \emph{economic}
    discretizations which lead to efficient numerical solutions. To the best of
    our knowlegde, only few key contributions on this topic have been made in
    recent years and there remains potential for improvement
    towards~\ref{itm:cdr:3}. In many applications only certain quantities are of
    interest. Adaptation with respect to these quantities improves efficiency
    significantly. Adjoint-based error estimation provides the tools for this
    goal-oriented approach. Specifically, we use the Dual Weighted Residual (DWR) method~\cite{BeRa01}.
    DWR has been applied successfully to parabolic convection-dominated problems
    in previous work by some of the
    authors~\cite{bruchhauser2023cost,bause2021flexible,BruchhaeuserB22}.
    Goal oriented error estimation with isotropic mesh adaptivity for space time problems using the DWR method can be found in \cite{BruchhaeuserSV08,BruchhaeuserBR12,endtmayer_goal-oriented_2024,RoThiKoeWi23,endtmayertouch2024goal,ThiWi24,EnLaRiSchaWi2024}.
    A framework for mesh adaptation with a high order DG discretization has been
    presented~\cite{yano_optimization-based_2012}.
    In~\cite{caplan_four-dimensional_2020} and references therein, these
    concepts have been generalized to anisotropic space-time adaptation.
    Furthermore, adaptively refined meshes are essential for numerical
    stability: On globally refined meshes even
    stabilization schemes fail to reduce oscillations near steep gradients
    (cf.~\cite{BruchhaeuserJS08}).

    The computational efficiency of AFEM can be improved through
    \emph{anisotropic} refinement along dominant error directions. This approach
    has been successfully applied to many fields,
    e.\,g.~\cite{chakraborty_anisotropic_2024,zander_anisotropic_2022,caplan_four-dimensional_2020,richter_anisotropic_2013,yano_optimization-based_2012,richter_posteriori_2010,leicht_error_2010,formaggia_anisotropic_2004,formaggia_anisotropic_2001,formaggia2001new,picasso2003anisotropic,IoanToulopSpaceTime-SUPGAnisotropic_2024}.
    Anisotropic mesh refinement applied to CDR problems can be found
    in~\cite{knobloch2025adaptive-aniso,yano_optimization-based_2012,georgoulis_discontinuous_2007,formaggia_anisotropic_2001,apel_anisotropic_1996,kornhuber_adaptive_1990}.
    For a general overview of anisotropic finite elements we refer to the
    monograph of Apel~\cite{apel_anisotropic_1999}.

    To address~\ref{itm:cdr:3} and achieve \emph{economic} discretizations, we
    propose an \emph{anisotropic space-time goal-oriented mesh adaptation}
    framework. Motivated by ideas of
    Richter~\cite{richter_posteriori_2010,richter_anisotropic_2013}, we combine
    goal-oriented error estimation with anisotropic refinement. Our proposed
    method provides directional error indicators in space and time with respect
    to a user-defined goal. Thereby, the directional separation of elementwise
    indicators naturally leads to anisotropic mesh refinement.

    Our approach automatically identifies dominant error directions. This
    represents one of the key challenges to other approaches using anisotropic mesh refinement. To this end, a frequently used approach involves the analysis
    of the Hessian matrix~\cite{leicht_error_2010,formaggia_anisotropic_2004,formaggia_anisotropic_2001}.
    Its eigenvectors indicate the major refinement directions, while the
    eigenvalues quantify the strength of the anisotropy. This method has been
    successfully applied in various settings. However, Hessian-based approaches
    face two key limitations. First, it is unclear whether refinement should be
    based on the Hessian of the primal solution \( u \) or the adjoint solution
    \( z \), as anisotropic influences may be caused by both variables, as explained in~\cite{richter_posteriori_2010}. Second, the Hessian is most effective
    when using linear finite elements. For higher-order finite elements, the
    Hessian does not provide reliable error information since second-order
    curvature is already captured by the method~\cite[Sec.
    8.4]{BruchhaeuserR17}. By directly estimating directional errors within the
    DWR framework these limitations are circumvented and a refinement
    strategy that accounts for both primal and adjoint residuals is achieved.

    The contribution of this work is twofold: (1) we develop a framework for
    anisotropic space-time mesh adaptation which captures directional
    solution features, and (2) we integrate goal-oriented error estimation to
    refine the mesh where it most impacts the goal, leading to computationally
    efficient and economic space-time discretizations. The implementation builds
    on the approaches presented
    in~\cite{KoeBruBau2019,bause2021flexible,bruchhauser2023cost}. Notably,
    aside from a few extra interpolation steps, estimating directional errors
    incurs no additional computational cost compared to the DWR method used in~\cite{KoeBruBau2019,bause2021flexible,bruchhauser2023cost}.

    The remainder of this work is organized as follows: In
    Sec.~\ref{sec:model-variational}, we introduce the CDR equation and its
    variational formulation. Sec.~\ref{sec:st-disc} outlines the space-time
    discretization and reviews the classical isotropic error representation as a
    starting point for further development. The main contribution is presented in
    Sec.~\ref{sec:anisotropic-method}, where we detail an anisotropic error
    estimation technique. In Sec.~\ref{sec:algorithm}, we describe the algorithm
    for goal-oriented anisotropic mesh refinement. Finally,
    Sec.~\ref{sec:8:numerical_examples} illustrates and validates our approach
    by numerical examples, including comparisons with isotropic refinement.

\section{Model Problem and Variational Formulation}\label{sec:model-variational}
As a model problem, we study the following time-dependent
CDR equation:
\begin{equation}
\label{eq:1:CDRoriginal}
\begin{array}{r@{\;}c@{\;}l@{\;} @{\,\,}l @{\,\,}l @{\,}l}
\partial_t u
- \nabla \cdot  (\varepsilon \nabla u)
+ \boldsymbol{b} \cdot \nabla u
+ \alpha u & = & f & \mbox{in } & \mathcal{Q} & = \Omega \times I \,,
\\
 u & = & u_D & \mbox{on } & \Sigma_D & = \Gamma_D \times I\,,
\\
 \varepsilon \nabla u \cdot \boldsymbol{n} & = & u_N & \mbox{on } & \Sigma_N & = \Gamma_N \times I\,,
\\
u(0) & = & u_{0} & \mbox{on } & \Sigma_0 & = \Omega\times \{0\} \,,
\end{array}
\end{equation}
in the space-time domain $\mathcal{Q}$, where $\Omega\subset \mathbb{R}^d$,
with $d=2$ or $d=3$, is a polygonal or polyhedral bounded domain with Lipschitz
boundary $\partial\Omega$ and $I=(0,T), 0 < T < \infty$, is a finite
time interval.
Here, $\partial\Omega = \Gamma_D \cup \Gamma_N\,,\Gamma_D \neq
\emptyset$ denotes the partition of the boundary with outer unit normal vector
$\boldsymbol{n}$, where $\Gamma_D$ denotes the
Dirichlet part and $\Gamma_N$ the Neumann part, respectively.
Furthermore, let $V:=\big\{v \in H^1(\Omega)|v_{|\Gamma_D}=0\big\}$ and $V'$ denotes the adjoint space of $V$.
To ensure the well-posedness of Eq.~\eqref{eq:1:CDRoriginal} we
assume that $0 < \varepsilon \leq 1$ is a constant diffusion coefficient,
$\boldsymbol{b} \in L^{\infty}\big(I;W^{1,\infty}(\Omega)^d\big)$
is the flow field or convection field,
$\alpha \in L^{\infty}\big(I;L^{\infty}(\Omega)\big)$
is a non-negative ($\alpha \geq 0$) reaction coefficient,
$u_0\in L^2(\Omega)$
is a given initial condition,
$f \in L^{2}(I;V')$ is a given source of the unknown scalar quantity $u$,
$u_D \in L^{2}(I;H^{\frac{1}{2}}(\Gamma_D))$ is a given function specifying the
Dirichlet boundary condition, and $u_N \in L^{2}(I;H^{-\frac{1}{2}}(\Gamma_N))$ is a given
function specifying the Neumann boundary condition.
Furthermore, it will be assumed that either $\nabla \cdot \boldsymbol{b}(\boldsymbol{x},t) = 0$
and $\alpha(\boldsymbol{x},t) \geq 0$, or there exists a
positive constant $c_0$ such that
$
\alpha(\boldsymbol{x},t)-\frac{1}{2} \textnormal{div}\;\boldsymbol{b}(\boldsymbol{x},t) \geq c_0 > 0
\;\;\forall (\boldsymbol{x},t) \in \bar{\Omega}\times \bar{I}\,,
$
which are standard assumptions for convection-dominated equations of type~\eqref{eq:1:CDRoriginal}, cf., e.g.,~\cite{BruchhaeuserAJ15,BruchhaeuserRST08}.

Henceforth, for the sake of simplicity, we deal with homogeneous Dirichlet
boundary values $u_D = 0$ on $\Gamma_D=\partial\Omega$ only. This implies that here $V\coloneqq H_0^1(\Omega)$. In the numerical examples in Sec.~\ref{sec:8:numerical_examples}, we also
consider more general boundary conditions, which are incorporated as described in~\cite[Ch.~3.3]{BaRa03}.

It is well known that problem \eqref{eq:1:CDRoriginal} along with the above
conditions admits a unique weak solution
$u \in \mathcal{V}\coloneqq \big\{v\in L^{2}\big(I; V\big)\big|
\;\partial_{t}v\in L^{2}(I;V')\big\}\,,$ that satisfies the following
variational formulation; cf., e.g.~\cite{BruchhaeuserRST08,BruchhaeuserJKN18}.
\begin{Problem}

  Find $u \in \mathcal{V}$ such that
  \begin{equation}
    \label{eq:2:WeakCDRsteady}
    A(u)(\varphi) = F(\varphi) \quad \forall \varphi \in \mathcal{V}\,,
  \end{equation}
  where the bilinear form $A:\mathcal{V}\times \mathcal{V} \rightarrow \mathbb{R}$
  and the linear form $F:\mathcal{V}\rightarrow \mathbb{R}$ are defined by
  \begin{eqnarray}
    \label{eq:3:BilinearformA}
    A(u)(\varphi) & \coloneqq & \int_{I}\big\{(\partial_t u,\varphi)
                                +a(u)(\varphi)
                                \big\} \mathrm{d} t
                                + (u(0),\varphi(0))\,,
    \\[1ex]
    \label{eq:4:LinearformF}
    F(\varphi) & \coloneqq & \int_I(f,\varphi)\;\mathrm{d}t
                             +(u_0,\varphi(0))\,,
  \end{eqnarray}
  with the inner bilinear form $a:V \times V \rightarrow \mathbb{R}$,
  given by
  \begin{equation}
    \label{eq:5:aBilinearform}
    a(u)(\varphi)\coloneqq (\varepsilon\nabla u, \nabla \varphi)
    +(\boldsymbol{b}\cdot \nabla u, \varphi) + (\alpha u,\varphi)\,.
  \end{equation}
\end{Problem}
We note that the initial condition is incorporated into the variational problem. The weak formulation, given by Eq.~\eqref{eq:2:WeakCDRsteady}, is now the
starting point for the variational discretization in space and time using Galerkin
finite element methods.

Later the aim of this work is the accurate approximation of a specific quantity of
interest \( J(u) \), defined by a functional \( J \in
\mathcal{C}^3(\mathcal{V},\mathbb{R}) \), where \( u \in \mathcal{V} \) solves
the model problem~\eqref{eq:2:WeakCDRsteady}.

\section{Space-Time Discretization and Isotropic Error Representation}\label{sec:st-disc}
In this section, we present the space-time finite element discretization of
our model problem including SUPG stabilization and, for completeness, review a goal-oriented error representation
based on the DWR method for the isotropic case.
\subsection{Discretization in Time}\label{sec:3.2:disc_time}
For the discretization in time we use a discontinuous Galerkin method $dG(r)$
with an arbitrary polynomial degree $r\geq0$.
Let $\mathcal T_{\tau}$ be a partition of the closure of the time domain
$\bar{I}=[0,T]$ into left-open subintervals $I_n\coloneqq (t_{n-1},t_n]$, $n=1,\dots,N$,
with $0=:t_0<t_1<\dots<t_N\coloneqq T$ and time step sizes $\tau_n=t_n-t_{n-1}$ and the global time discretization
parameter $\tau=\max_{n}\,\tau_{n}$.
Since the set of time intervals \( I_n \) is finite, it is natural to partition the global space-time cylinder \( \mathcal{Q} = \Omega \times I \) into space-time slabs defined as \( \mathcal{Q}_n = \Omega \times I_n \).
On the subintervals $I_n$, we define the time-discrete function space
$\mathcal{V}_{\tau}^{r}$
\begin{equation}
\label{eq:6:Def_V_tau_dGr}
 \begin{aligned}
\mathcal{V}_{\tau}^{r} \coloneqq
 \Big\{u_{\tau}\in L^{2}(I; V)\big|
 u_{\tau}|_{I_{n}}\in \mathcal{P}_{r}(I_{n}; V),
 u_{\tau}(0)\in L^2(\Omega), n=1,\dots,N
\Big\}\,,
\end{aligned}
\end{equation}
where $\mathcal{P}_{r}(\bar{I}_{n}; V)$ denotes the space of all
polynomials in time up to degree $r\geq0$ on $I_n$ with values in $V\,.$
For some discontinuous in time function $u_{\tau}\in \mathcal{V}_{\tau}^{r}$
we define the limits $u_{\tau,n}^{\pm}$ from above and below of $u_{\tau}$ at
$t_n$ as well as their jump at $t_n$ by
\begin{displaymath}
\begin{array}{lcrclcr}
u_{\tau,n}^{\pm}
& \coloneqq &
\displaystyle\lim_{t\to t_n\pm0} u_\tau(t) \,,
&
[u_{\tau}]_{n} & \coloneqq & u_{\tau,n}^{+}
-u_{\tau,n}^{-} \,.
\end{array}
\end{displaymath}
Then, the semi-discrete in time scheme of Eq.~\eqref{eq:2:WeakCDRsteady} reads
as follows:
\textit{
Find $u_{\tau} \in \mathcal{V}_{\tau}^{r}$ such that
}
\begin{equation}
\label{eq:7:dGDiscTime}
A_{\tau}(u_{\tau})(\varphi_{\tau}) =
F_{\tau}(\varphi_{\tau}) \quad \forall \varphi_{\tau}
\in \mathcal{V}_{\tau}^{r}\,,
\end{equation}
\textit{
where the bilinear form $A_{\tau}(\cdot)(\cdot)$ and the linear form
$F_{\tau}(\cdot)$ are defined by
}
\begin{eqnarray}
\label{eq:8:BilinearFormAtau}
A_{\tau}(u_{\tau})(\varphi_{\tau}) & \coloneqq &
\displaystyle\sum_{n=1}^{N}\int_{I_n}
\big\{(\partial_t u_{\tau},\varphi_{\tau})
+a(u_{\tau})(\varphi_{\tau})
\big\} \mathrm{d} t
+ \displaystyle\sum_{n=2}^N
([u_{\tau}]_{n-1},\varphi_{\tau,n-1}^+ )
+ (u_{\tau,0}^{+},\varphi_{\tau,0}^{+})\,,
\\[1ex]
\label{eq:9:LinearFormFtau}
F_{\tau}(\varphi_{\tau}) & \coloneqq & \displaystyle\int_I(f,\varphi)\;\mathrm{d}t
+(u_0,\varphi_{\tau,0}^{+})\,,
\end{eqnarray}
with the inner bilinear form $a(\cdot)(\cdot)$ being defined by
Eq.~\eqref{eq:5:aBilinearform}.

\subsection{Discretization in Space and Stabilization}
\label{sec:3.3:disc_space}
Next, we describe the Galerkin finite element approximation in space for the
semi-discrete time scheme~\eqref{eq:7:dGDiscTime}. We use Lagrange-type finite
element spaces of continuous, piecewise polynomial functions. The spatial
discretization is based on a decomposition \( \mathcal{T}_h \) of the domain
\( \Omega \) into disjoint elements \( K \), such that
\( \overline{\Omega} = \cup_{K\in\mathcal{T}_h} \overline{K} \). For
\( d=2,3 \), we use quadrilateral and hexahedral elements, respectively. Each
element \( K \in \mathcal{T}_h \) is mapped from the reference element
$\hat{K}=(0,1)^d$ via an
iso-parametric transformation \( \boldsymbol T_K: \hat{K} \to K \) satisfying
\( \det(\boldsymbol T_K)(\hat{x}) > 0 \) for all \( \hat{x} \in (0,1)^d \).
Following~\cite{richter_posteriori_2010}, we decompose
\begin{equation}\label{eq:TK}
  \boldsymbol T_K \coloneqq \boldsymbol R_K \circ\boldsymbol S_{c,K} \circ\boldsymbol S_{h,K} \circ\boldsymbol P_K,
\end{equation}
where \(\boldsymbol R_K \) is a rotation and translation, \(\boldsymbol S_{c,K} \) an anisotropic
scaling, \(\boldsymbol S_{h,K} \) a shearing, and \(\boldsymbol P_K \) is a nonlinear component. To
account for anisotropic elements, we relax the standard shape-regularity
assumptions and require only uniform boundedness of \(\boldsymbol S_{h,K} \) and \(\boldsymbol P_K \)
for all \( K \in \mathcal{T}_h \). The element diameter is denoted by \( h_K \),
with the global discretization parameter defined as
\( h \coloneqq \max_{K\in\mathcal{T}_h} h_K \).

Our mesh adaptation procedure employs local refinement with hanging nodes.
Global conformity of the finite element approach is preserved by eliminating
degrees of freedom at hanging nodes via interpolation between neighboring
regular nodes; see~\cite[Ch.~4.2]{BaRa03} and~\cite{BruchhaeuserCO84}.

\begin{Definition}
  On a subsect $\mathcal Z_h\subseteq \mathcal T_h$, we define the discrete finite element space
  $$
  V_h^{p}(\mathcal Z_h)\coloneqq
  \big\{v\in C(\overline{\Omega})\mid v_{|K}
  \in \Q{p}(K)\,,\forall K\in\mathcal Z_h
  \big\}\cap V\,,
  $$
  where, $\Q{p}(K)$ is the mapped finite element from
  \[
    \hat{Q}_{p}(\hat{K}) \coloneqq \bigotimes_{i=1}^{d} \mathcal{P}_{p}([0,1])\,,
  \]
  by~\eqref{eq:TK}, where \(\mathcal{P}_{p}([0,1])\) is the space of univariate
  polynomials of degree $p$. On $\mathcal{T}_h$ we put
  $$V_h^{p}=V_h^{p}(\mathcal{T}_h)\,.$$
\end{Definition}

Let the fully discrete function space be given by
\begin{equation}
\label{eq:10:Def_V_tau_h_dGr_p}
\mathcal{V}_{\tau h}^{r,p} \coloneqq \Big\{
u_{\tau h}\in X_{\tau}^{r} \big|
u_{\tau h}|_{I_n} \in \mathcal{P}_r(I_n;V_h^{p})
\,, u_{\tau h}(0) \in V_h^{p},
n=1,\dots,N
\Big\}
\subseteq L^{2}(I; V)\,.
\end{equation}
We note that the spatial finite element space $V_h^{p}$ is allowed to be
different on all subintervals $I_n$ which is natural in the context of a
discontinuous Galerkin approximation of the time variable and allows dynamic
mesh changes in time.
Due to the conformity of $V_h^{p}$ we get $\mathcal{V}_{\tau h}^{r,p}\subseteq
\mathcal{V}_{\tau}^{r}$.
Now, the fully discrete discontinuous in time scheme reads as follows:
\textit{Find $u_{\tau h} \in \mathcal{V}_{\tau h}^{r,p}$ such that}
\begin{equation}
\label{eq:11:A_tau_h_u_phi_eq_F_phi}
 A_{\tau}(u_{\tau h})(\varphi_{\tau h})
=
F_\tau(\varphi_{\tau h})
\quad \forall \varphi_{\tau h} \in \mathcal{V}_{\tau h}^{r,p}\,,
\end{equation}
\textit{with} $A_{\tau}(\cdot)(\cdot), a(\cdot)(\cdot)$ \textit{and} $F_\tau(\cdot)$
\textit{being defined in \eqref{eq:8:BilinearFormAtau},\eqref{eq:5:aBilinearform} and
\eqref{eq:9:LinearFormFtau}, respectively.}

In this work, we focus on convection-dominated problems with small diffusion
parameter $\varepsilon >0$. Then, the finite element approximation needs to be
stabilized in order to reduce spurious and non-physical oscillations of the
discrete solution arising close to layers.
Here, we apply the streamline upwind Petrov-Galerkin (SUPG) method
\cite{BruchhaeuserHB79,BruchhaeuserBH82}.
Existing convergence analyses in the natural norm of the underlying scheme
including local and global error bounds can be found, for instance, in~\cite[Ch. III.4.3]{BruchhaeuserRST08}.
The stabilized variant of the fully discrete discontinuous in time scheme
then reads as follows:
\textit{
Find $u_{\tau h} \in \mathcal{V}_{\tau h}^{r,p}$ such that
}
\begin{equation}
\label{eq:12:StabDGFully}
A_{S}(u_{\tau h})(\varphi_{\tau h}) =
F_{\tau}(\varphi_{\tau h}) \quad \forall \varphi_{\tau h}
\in \mathcal{V}_{\tau h}^{r,p}\,,
\end{equation}
\textit{
where the linear form $F_{\tau}(\cdot)$ is defined by
Eq.~\eqref{eq:9:LinearFormFtau} and the stabilized bilinear form
$A_{S}(\cdot)(\cdot)$ is given by}
\begin{equation}
\label{eq:13:ASutauhDG}
A_{S}(u_{\tau h})(\varphi_{\tau h}) \coloneqq
A_{\tau}(u_{\tau h})(\varphi_{\tau h})
+ S(u_{\tau h})(\varphi_{\tau h})\,.
\end{equation}
\textit{
Here, the SUPG stabilization bilinear form
}
$S(\cdot)(\cdot)$
\textit{
is defined by
}
\begin{equation}
\label{eq:14:StabDGSutauh}
\begin{array}{r@{\;}c@{\;}l@{\;}}
S(u_{\tau h})(\varphi_{\tau h}) & \coloneqq &
\displaystyle
\sum_{n=1}^{N}\int_{I_n}\sum_{K \in \mathcal{T}_{h, n}}
\delta_K\big(
r(u_{\tau h}),
\boldsymbol{b} \cdot \nabla \varphi_{\tau h}\big)_K \,\mathrm{d} t
\\[3ex]
& &
\displaystyle
+ \sum_{n=2}^{N}\sum\limits_{K\in\mathcal{T}_{h, n}}
\delta_K
\big(\left[u_{\tau h}\right]_{n-1},
\boldsymbol{b} \cdot \nabla \varphi_{\tau h,n-1}^+\big)_{K}
\\[3ex]
& &
+ \displaystyle\sum\limits_{K\in\mathcal{T}_{h, 1}}
\delta_K \big(u_{\tau h,0}^+ - u_0,
\boldsymbol{b} \cdot \nabla \varphi_{\tau h,0}^{+} \big)_{K}\,,
\end{array}
\end{equation}
\textit{
where the residual term $r(\cdot)$
is defined by}
\begin{equation}
\label{eq:15:rutauh}
r(u_{\tau h}) \coloneqq \partial_{t} u_{\tau h}
- \nabla\cdot\left(\varepsilon\nabla u_{\tau h}\right)
+ \boldsymbol{b} \cdot \nabla u_{\tau h}
+ \alpha u_{\tau h}
- f \,.
\end{equation}
\begin{Remark}
\label{rem:2:StabParameter}
The proper choice of the stabilization parameter $\delta_K$ is an important
issue in the application of the SUPG approach; cf., e.g.,
\cite{BruchhaeuserJN11,BruchhaeuserJS08,BruchhaeuserJKN18}
and the discussion therein. For time-dependent CDR
problems an optimal error estimate for $\delta_K=\mathrm{O}(h)$ is derived
in~\cite{BruchhaeuserJN11}. In this work, we chose $\delta_{K}=\delta_0 h_K$, where $\delta_0=0.1$ and $h_K=\sqrt[d]{|K|}$ is the cell diameter.
\end{Remark}

%
\subsection{Isotropic Goal Oriented Error Estimation}
\label{sec:3.4:IsotropicEE}
As the fundamental ideas translate to anisotropic error representation, we present an isotropic a posteriori error representation for the stabilized problem~\eqref{eq:12:StabDGFully} based on the Dual Weighted Residual (DWR)
method~\cite{BeRa01,BaRa03,BruchhaeuserBR96}.
%
This representation is given in terms of a user-chosen goal functional
\begin{equation*}
  J(u)=\int_0^T J_1(u(t))\mathrm{d}t
  + J_2(u(T))\,,
\end{equation*}
where $J_1$ and $J_2$ are three times differentiable functionals and each of
them may be zero; cf.~\cite{BruchhaeuserSV08,BruchhaeuserBR12}. Since we focus
on anisotropic error estimation, we restrict this section to the main result and
refer to~\cite[Ch.~4]{BruchhaeuserB22} for a detailed derivation in the context
of CDR equations.

For the error representations in Theorem.~\ref{Thm:3.3}, we introduce the following Lagrangian
functionals
$\mathcal{L}: \mathcal{V}\times \mathcal{V} \rightarrow \mathbb{R}$,
$\mathcal{L}_\tau: \mathcal{V}_{\tau}^{r} \times \mathcal{V}_{\tau}^{r}
\rightarrow \mathbb{R}$, and
$\mathcal{L}_{\tau h}:
\mathcal{V}_{\tau h}^{r,p} \times \mathcal{V}_{\tau h}^{r,p}
\rightarrow \mathbb{R}$ by
\begin{subequations}
  \label{eq:3:3:Def_L_u_z_Def_L_tau_u_z_Def_L_tau_h_u_z}
  \begin{align}
    \label{eq:3:3:Def_L_u_z}
    \mathcal{L}(u,z) & \coloneqq  J(u)
                       + F(z)
                       - A(u)(z)\,,
    \\
    \label{eq:3:3:Def_L_tau_u_z}
    \mathcal{L}_{\tau}(u_\tau,z_\tau) & \coloneqq
                                        J(u_{\tau}) + F_\tau(z_{\tau})
                                        - A_{\tau}(u_{\tau})(z_{\tau})\,,
    \\
    \label{eq:3:3:Def_L_tau_h_u_z}
    \mathcal{L}_{\tau h}(u_{\tau h},z_{\tau h}) & \coloneqq
                                                  J(u_{\tau h})
                                                  + F_\tau (z_{\tau h})
                                                  - A_S(u_{\tau h})(z_{\tau h})=\mathcal{L}_{\tau}(u_{\tau h}),z_{\tau h}))-S(u_{\tau h}))(z_{\tau h}))\,.
  \end{align}
\end{subequations}
In~\eqref{eq:3:3:Def_L_u_z_Def_L_tau_u_z_Def_L_tau_h_u_z}, the Lagrange
multipliers $z$, $z_\tau$, and $z_{\tau h}$ are called adjoint variables in
contrast to the primal variables $u$, $u_\tau$, and $u_{\tau h}$;
cf.~\cite{BruchhaeuserBR12,BeRa01}.
The directional derivatives of the Lagrangian functionals
(G\^{a}teaux derivatives) with respect to their second argument yields the
primal
problems~\eqref{eq:2:WeakCDRsteady},~\eqref{eq:7:dGDiscTime},~\eqref{eq:12:StabDGFully},
while the directional derivative with respect to their first argument yields the
adjoint problems given by
\begin{subequations}
\label{eq:3:9:A_tau_prime_u_phi_z_eq_J_prime_u_phi_A_S_prime_u_phi_z_eq_J_prime_u_phi}
\begin{align}
  \label{eq:3:6:A_prime_u_phi_z_eq_J_prime_u_phi}
  A^{\prime}(u)(\varphi,z)
  &=
  J^{\prime}(u)(\varphi)
  \quad \forall \varphi \in \mathcal{V}\,,
  \\
  \label{eq:3:9:A_tau_prime_u_phi_z_eq_J_prime_u_phi}
  A_{\tau}^{\prime}(u_{\tau})(\varphi_{\tau},z_{\tau})
  & =
    J^{\prime}(u_{\tau})(\varphi_{\tau})
    \quad \forall \varphi_{\tau}\in \mathcal{V}_{\tau}^{r}\,,
  \\
  \label{eq:3:9:A_S_prime_u_phi_z_eq_J_prime_u_phi}
  A_{S}^{\prime}(u_{\tau h})(\varphi_{\tau h},z_{\tau h})
  & =
    J^{\prime}(u_{\tau h})(\varphi_{\tau h})
    \quad \forall \varphi_{\tau h}\in \mathcal{V}_{\tau h}^{r,p}\,,
\end{align}
\end{subequations}
where the definitions of $A^{\prime}(\cdot)(\cdot,\cdot)$, $A_{\tau}^{\prime}(\cdot)(\cdot, \cdot)$,
$A_{S}^{\prime}(\cdot)(\cdot, \cdot)$ are given in the appendix.
\begin{theorem}
\label{Thm:3.3}
Let $\{u,\,z\}\in \mathcal{V} \times \mathcal{V}$,
$\{u_{\tau},z_{\tau}\}
\in
\mathcal{V}_{\tau}^{r} \times \mathcal{V}_{\tau}^{r}$,
and
$\{u_{\tau h},z_{\tau h}\}
\in \mathcal{V}_{\tau h}^{r,p} \times \mathcal{V}_{\tau h}^{r,p}$
denote the stationary points of
$\mathcal{L}, \mathcal{L}_{\tau}$, and $\mathcal{L}_{\tau h}$
on the different levels of discretization, i.e.,
\begin{displaymath}
\begin{aligned}
\mathcal{L}^{\prime}(u,z)(\delta u, \delta z)
= \mathcal{L}_{\tau}^{\prime}(u,z)(\delta u, \delta z)
& = 0 \quad
\forall \{\delta u,\delta z\}\in \mathcal{V} \times \mathcal{V}\,,
\\
\mathcal{L}_{\tau}^{\prime}(u_{\tau},z_{\tau})
(\delta u_{\tau}, \delta z_{\tau})
& = 0
\quad \forall \{\delta u_{\tau},\delta z_{\tau}\}
\in \mathcal{V}_{\tau}^{r} \times \mathcal{V}_{\tau}^{r}\,,
\\
\mathcal{L}_{\tau h}^{\prime}(u_{\tau h},z_{\tau h})
(\delta u_{\tau h}, \delta z_{\tau h})
& = 0
\quad \forall \{\delta u_{\tau h},\delta z_{\tau h}\}
\in \mathcal{V}_{\tau h}^{r,p} \times \mathcal{V}_{\tau h}^{r,p}\,.
\end{aligned}
\end{displaymath}
Then, for the discretization errors in space and time we get
the representation formulas
\begin{subequations}
\label{eq:3:16}
\begin{align}
\label{eq:3:16a:J_u_minus_J_u_tau}
J(u)-J(u_{\tau}) & =
\frac{1}{2}\rho_{\tau}(u_{\tau})(z-\tilde{z}_{\tau})
+ \frac{1}{2}\rho_{\tau}^{\ast}(u_{\tau},z_{\tau})
(u-\tilde{u}_{\tau})
+ \mathcal{R}_{\tau}\,,
\\
\label{eq:3:16b:J_u_tau_minus_J_u_tau_h}
J(u_{\tau})-J(u_{\tau h}) & =
\frac{1}{2}\rho_{\tau}(u_{\tau h})(z_{\tau}-\tilde{z}_{\tau h})
+ \frac{1}{2}
\rho_{\tau}^{\ast}(u_{\tau h},z_{\tau h})
(u_{\tau}-\tilde{u}_{\tau h})
\\
\nonumber
& \qquad
+ \frac{1}{2} S(u_{\tau h})(\tilde{z}_{\tau h}+z_{\tau h})
+ \frac{1}{2} S^{\prime}(u_{\tau h})(\tilde{u}_{\tau h}-u_{\tau h},z_{\tau h})
+ \mathcal{R}_{h}\,,
\end{align}
\end{subequations}
where $\rho_{\tau}$ and $\rho_{\tau}^{\ast}$ are the primal and adjoint residuals based on the semi-discrete
in time schemes, respectively, given by
\begin{subequations}
  \begin{align}
    \label{eq:3:15:primal_dual_residuals}
    \rho_{\tau}(u)(\varphi)  &\coloneqq
    \mathcal{L}_{\tau,z}^{\prime}(u,z)(\varphi)=F_{\tau}(\varphi) - A_{\tau}(u)(\varphi)\,,\\
    \rho_{\tau}^{\ast}(u,z)(\varphi)
    &\coloneqq
    \mathcal{L}_{\tau,u}^{\prime}(u,z)(\varphi)= J^{\prime}(u)(\varphi) - A^{\prime}_{\tau}(u)(\varphi,z)
    \,.
  \end{align}
\end{subequations}
Here,
$\{\tilde{u}_{\tau},\tilde{z}_{\tau}\}\in \mathcal{V}_{\tau}^{r}
\times \mathcal{V}_{\tau}^{r}$,
and
$\{\tilde{u}_{\tau h},\tilde{z}_{\tau h}\} \in
\mathcal{V}_{\tau h}^{r,p} \times \mathcal{V}_{\tau h}^{r,p}$
can be chosen arbitrarily and $\mathcal{R}_{\tau}, \mathcal{R}_{h}$
are higher-order remainder terms with respect to the errors
$u-u_{\tau}, z-z_{\tau}$ and $u_{\tau}-u_{\tau h}, z_{\tau}-z_{\tau h}$, respectively.
\end{theorem}
\subsection{Isotropic Patch-Wise Higher Order Interpolation} \label{subsection: IsoHighInt}
In this section we introduce the isotropic patch-wise higher order
interpolation. This interpolation technique is widely used in goal-oriented
error estimation (cf., e.\,g.,~\cite{BeRa01,BaRa03,RanVi2013,RiWi15_dwr,endtmayer2021reliability}).
We assume that $\mathcal{T}_h$ has been once refined, i.\,e.\
$\mathcal{T}_h=R(\mathcal{T}_{2h})$, where $R$ denotes the uniform refinement of
a mesh.

\begin{figure}[htb]
  \begin{minipage}{0.49\textwidth}\centering
    {\input{Tikz/FigureQ1.tex} }
  \end{minipage} \hfill%
  \begin{minipage}{0.49\textwidth}\centering
    {\input{Tikz/FigureQ2.tex} }
  \end{minipage}
  \caption{DoFs of $\Q{1}$ finite elements on the patch $K_{2h}^P$ (left) and a
    the DoFs of a $\Q{2}$ finite element on the element $K_{2h}$ (right).\label{Fig: Patch QIso Nodes}}
\end{figure}
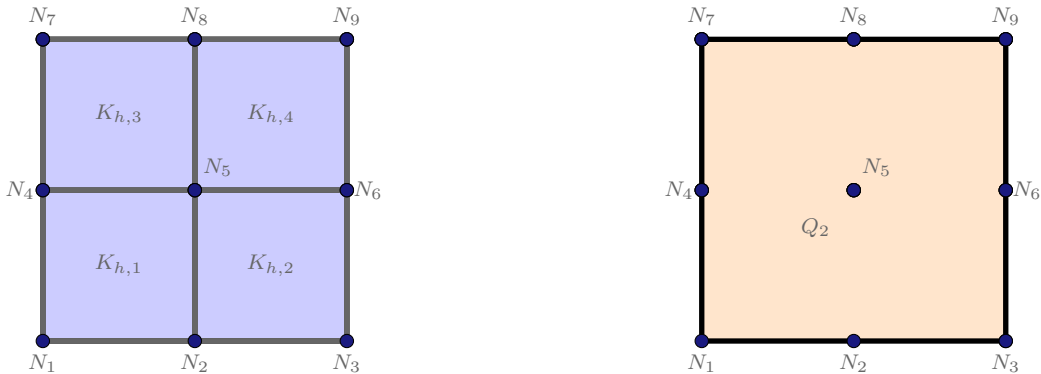
To an element $K_{2h} \in \mathcal{T}_{2h}$, we uniquely associate the patch
\begin{align}
  K_{2h}^P&\coloneqq \{K \in \mathcal{T}_h: |K \cap K_{2h}|\not=0\}\,,
  \intertext{and define the unique partitioning of $\mathcal{T}_h$ into patches}
  \mathcal T_{2h}^P&\coloneqq \bigcup_{K_{2h}\in \mathcal{T}_{2h}}K_{2h}^P\,.
\end{align}
We note that, the degrees of freedom (DoFs) of $V_h^p(K_{2h}^P)$
coincide with the DoFs in $V_h^{2p}(K_{2h})$. We therefore enumerate the DoFs on $V_h^p(K_{2h}^P)$ and $V_h^{2p}(K_{2h})$ identically, from $1$ to $N_{\text{DoF}}^P \coloneqq \dim V_h^p(K_{2h}^P)$.
For $\Q{1}$ and $\Q{2}$ finite elements this is viusalized in
Figure~\ref{Fig: Patch QIso Nodes}.
\begin{figure}[htb]
  \begin{minipage}[t]{0.49\textwidth}\centering
    \scalebox{0.75}{\input{Tikz/PatchQ11.tex} }
    \subcaption{The function $v_h$. \label{Fig: PatchQ11 Function}}
  \end{minipage} \hfill%
  \begin{minipage}[t]{0.49\textwidth}\centering
    \scalebox{0.75}{\input{Tikz/PatchQ22HighInter.tex} }
    \subcaption{The resulting function $\mathrm{I}_{2h}^{\left(2p\right)}v_h$. \label{Fig: PatchQ22 Function}}
  \end{minipage}
  \caption{The action of $\mathrm{I}_{2h}^{\left(2p,\right)}$ of a $Q_1$ finite element function $v_h$ on the patch $K_{2h}^P$. \label{fig: action of IsoInt}}
\end{figure}
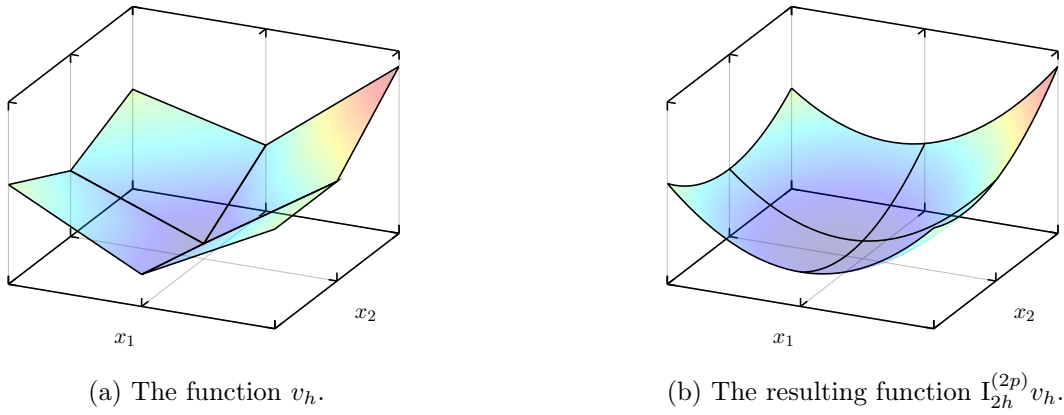

On the patches, we define the isotropic patch-wise higher order interpolation
$\mathrm{I}_{2h}^{(2p)}$ as in~\cite{BaRa03}.
\begin{Definition} \label{def: IsoHighInt}
  We define the bijective mapping $\mathrm{I}_{2h}^{\left(2p,K_{2h}\right)}:
  V_h^{p}(K_{2h}^P) \to V_h^{2p}(K_{2h})$ by
  \begin{equation*}
  \label{eq:def:I2h}
  \left(\mathrm{I}_{2h}^{\left(2p,K_{2h}\right)}v_h\right)(N_j)=v_h(N_j) \qquad
  \text{for } j=1,\dots,\,N_{\text{DoF}}^P\,,
\end{equation*}
for all DoFs $N_j\in K_{2h}^P$ and $v_h \in V_h^p(K_{2h}^P)$.
The global mapping $\mathrm{I}_{2h}^{\left(2p\right)}\colon V_h^p(\mathcal
T_h)\to V_{h}^{2p}(\mathcal T_{2h}^P)$ is given by
\begin{equation}\label{eq:iso-high-int}
  \mathrm{I}_{2h}^{\left(2p\right)} v_h(x)\coloneqq
  \left(\mathrm{I}_{2h}^{\left(2p,K_{2h}\right)}v_h\right)(x)\qquad \text{for }
  x \in K\subset K_{2h}^P\,,\, K_{2h}^P\in \mathcal{T}_{2h}^P\,.
\end{equation}
\end{Definition}
For later computations, we also need the restriction operator on the reference
element $\hat{K}$.
\begin{Definition}\label{def: iso element restriction}
  Let $\hat{v}^{2p} \in \hat{Q}_{2p}(\hat{K})$. Then, the restriction
  $\hat R^{2p,p}:\hat{Q}_{2p}(\hat{K}) \to \hat{Q}_{p}(\hat{K})$ is defined as the
  unique function defined by the condition
  \begin{equation}\label{def: iso space restriction ref}
    \left(\hat R^{2p,p}\hat v_{2p}\right)(\hat N_j) = \hat v_{2p}(\hat N_j)\qquad \text{for } j=1,\dots,\,(p+1)^d\,,
  \end{equation}
  where $\hat N_j$ is a DoF on $\hat{Q}_{p}(\hat{K})$. We call
  $\hat R^{2p,p}$ the isotropic restriction operator on $\hat{K}$.
\end{Definition}

\begin{Definition}\label{def: iso space restriction}
For $v_h \in \mathcal{V}_{h}^{2p}$, we define the restriction operator
$R^p_h:\mathcal{V}_{h}^{2p} \to \mathcal{V}_{h}^{p}$ using the isoparametric
mapping~\eqref{eq:TK} and~\eqref{def: iso space restriction ref}, as
\begin{equation}\label{eq:iso-space-restriction}
  (R_h^p v_h)(x)\coloneqq \left(\hat R^{2p,p}\hat v_h\right)\left(\boldsymbol{T}_K^{-1}(x)\right), \quad x \in K\,,
\end{equation}
 with $\hat v_h =
\boldsymbol{T}_K^{-1}\circ v_h$.
\end{Definition}
\subsection{Local Error Indicators in Space and Time}
\label{sec:3.4.2:EI}
For the practical use of the error representations~\eqref{eq:3:16}
obtained in Theorem~\ref{Thm:3.3}, we need to replace unknown quantities by
computable ones. This requires multiple steps:
First, we neglect the remainder
terms  $\mathcal{R}_{\tau}$ and $\mathcal{R}_{h}$ in~\eqref{eq:3:16} as they are
of higher order.
Second, the unknown solutions in~\eqref{eq:3:16} are replaced by means of the approximated fully discrete solutions $u_{\tau h}\in\mathcal{V}_{\tau h}^{r,p}$
and $z_{\tau h}\in \mathcal{V}_{\tau h}^{r,q}$, with $q=2p$, respectively.
Third, the temporal and spatial weights have to be approximated in a suitable
way, see~\cite{BeRa01,BaRa03} for a short review of possible techniques.
Here, we choose the following approximations for the temporal and spatial weights
based on our previous numerical comparisons in~\cite{BruSchweBau18,BruchhaeuserBB24}:
\paragraph{Temporal Weights.} Approximate the temporal weights $u-\tilde{u}_{\tau}$ and $z-\tilde{z}_{\tau}$,
respectively, by means of a higher-order reconstruction using Gauss
quadrature points given by
\begin{displaymath}
u-\tilde{u}_{\tau} \approx \mathrm{E}_{\tau}^{(r+1)}u_{\tau h}-u_{\tau h}\,,
\qquad
z-\tilde{z}_{\tau} \approx \mathrm{E}_{\tau}^{(r+1)}z_{\tau h}-z_{\tau h}\,,
\end{displaymath}
using a reconstruction in time operator $\mathrm{E}_{\tau}^{(r+1)}$
thats lifts the solution to a piecewise
polynomial of degree ($r+1$) in time, cf.~\cite{BruchhaeuserBB24,BruchhaeuserB22} for further details.
\paragraph{Spatial Weights.} Using the patch-wise higher-order interpolation
operator~\eqref{eq:iso-high-int} the restriction operator~\eqref{eq:iso-space-restriction}, we approximate the spatial weights $u_{\tau}-\tilde{u}_{\tau h}$, $z_{\tau}-\tilde{z}_{\tau h}$, and $\tilde{z}_{\tau h}$ via patch-wise higher-order interpolation and higher-order finite element approach as
\begin{equation}
  u_{\tau}-\tilde{u}_{\tau h}  \approx  \mathrm{I}_{2h}^{(2p)}u_{\tau h}-u_{\tau h}\,,
  \qquad
  z_{\tau}-\tilde{z}_{\tau h}  \approx  z_{\tau h}-\mathrm{R}_{h}^{p}z_{\tau h}\,,
  \qquad
  \tilde{z}_{\tau h} \approx  \mathrm{R}_{h}^{p}z_{\tau h}\,.
\end{equation}

Taking the above into consideration, we then obtain the temporal and spatial error estimators
\begin{subequations}\label{eq: implemented error estimator}
  \begin{align}\label{eq: implemented error estimator-tau}
    \eta_{\tau}&\coloneqq\frac{1}{2}\rho_{\tau}(u_{\tau h})(\mathrm{E}_{\tau}^{(r+1)}z_{\tau h}-z_{\tau h})
    + \frac{1}{2}\rho_{\tau}^{\ast}(u_{\tau h},z_{\tau
    h})(\mathrm{E}_{\tau}^{(r+1)}u_{\tau h}-u_{\tau h})\,,\\
    \label{eq: implemented error estimator-h}
    \eta_{h}&\coloneqq\frac{1}{2}\rho_{\tau}(u_{\tau h})(z_{\tau h}-\mathrm{R}_{h}^{p}z_{\tau h})
    + \frac{1}{2}
    \rho_{\tau}^{\ast}(u_{\tau h},\mathrm{R}_{h}^{p}z_{\tau h})
              (\mathrm{I}_{2h}^{\left(2p\right)}u_{\tau h}-u_{\tau h})\, \nonumber
    \\
     &\qquad
    + \frac{1}{2} S(u_{\tau h})(z_{\tau h}+\mathrm{R}_{h}^{p}z_{\tau h})
      + \frac{1}{2} S^{\prime}(u_{\tau h})(\mathrm{I}_{2h}^{\left(2p\right)}u_{\tau h}-u_{\tau h},\mathrm{R}_{h}^{p}z_{\tau h})\,.
  \end{align}
\end{subequations}
Finally, the error indicators have to be represented in a localized form in
order to mark elements within the adaptive mesh refinement algorithm.
Here, we follow the approach of~\cite{BeRa01,BaRa03}, where we recast the
estimators in~\eqref{eq: implemented error estimator} as
\begin{equation}
  \label{eq:localized-errors-iso}
  \eta_{\tau}=\sum_{n=1}^{N}\sum_{K\in\mathcal{T}_h^n}\eta_{\tau}^{K,n},\qquad\eta_h=\sum_{n=1}^{N} \sum_{K\in\mathcal{T}_h^n} \eta_{h}^{K,n}\,.
\end{equation}
We refer to~\cite[Sec.~4.3.2]{BruchhaeuserB22} for a
detailed derivation of~\eqref{eq:localized-errors-iso} in the context of CDR equations.

\section{Anisotropic Method}\label{sec:anisotropic-method}
Convection-dominated problems often exhibit anisotropic solution features
(e.\,g.\ thin layers with steep gradients). Anisotropic mesh refinement is advantageous in these cases,
as interpolation error estimates indicate that approximation errors scale with
the mesh width in each coordinate direction multiplied by the corresponding
directional derivative norm (cf.~\cite{apel_anisotropic_1998}). Thus,
being able to identify dominant error directions a posteriori, offers the
potential to increase computational efficiency compared to isotropic refinement.
In this section, we split the isotropic error estimator $\eta_h$
(cf.~\eqref{eq: implemented error estimator-h}) into directional contributions such that
\begin{equation*}
  \eta_{h}\coloneqq\eta_{h,1}+ \eta_{h,2}+\eta_{h,\mathbb{E}}\,,
\end{equation*}
where $\eta_{h,\mathbb{E}}$ is an anisotropic remainder term.
Further, we describe the anisotropic mesh adaptation scheme based on the
directional error estimators.

\subsection{Anisotropic Finite Elements}
In order to split the error in directional contributions, we first
introduce finite elements with anisotropic polynomial degree.
Let \((p_1,\dots,p_d)\in\mathbb{N}^d\) be a multi-index.
On the reference element $\hat{K}$ the anisotropic polynomial space is defined as
\[
\hat{Q}_{p_1,\dots,p_d}(\hat{K}) \coloneqq \bigotimes_{i=1}^{d} \mathcal{P}_{p_i}([0,1]),
\]
where \(\mathcal{P}_{p_i}([0,1])\) denotes the space of univariate polynomials
of degree at most \(p_i\). The DoFs are associated with Gauss--Lobatto nodes.
The total number of DoFs per element is
\(
N_K^{p_1,\dots,p_d} = \prod_{i=1}^{d} (p_i+1).
\)
We note that this is a straightforward generalization of the isotropic case, where for
\(p_1=\cdots=p_d \coloneqq p\), we put $\hat{Q}_{p}$.
For \(i\in\{1,\dots,d\}\) and given integers \(p\) and \(q\), we define
\[
\hat{Q}^{p,q}_i \coloneqq \hat{Q}^{p_1,\dots,p_d}, \quad
p_j =
\begin{cases}
  p, & j\neq i \\
  q, & j=i
\end{cases}\,.
\]
Further we denote $N_i^{p,q}$ denotes the set of DoFs on the reference element
$\hat K$.
\subsection{Anisotropic Interpolation and Restriction Operations} \label{subsection: AnIsoHighInt}
In this subsection, we introduce anisotropic interpolation and restriction
operators, which are essential for computing the weights in the anisotropic
error estimator. Throughout the remainder of this section, we restrict ourselves
to the two-dimensional case, i.\,e. $\Omega \subset \mathbb{R}^2$.
\begin{Definition}
  For $p\leq q$, the restriction $\hat R_i^{q,p}\colon
  \hat{Q}_q(\hat{K}) \to \hat{Q}_i^{p,q}(\hat{K})$ is defined as the unique
  function such that, for $\hat{v}^q \in \hat{Q}_q(\hat{K})$,
  \begin{equation}
    \left(\hat R_i^{q,p}\hat v^q\right)(\hat N_j) = \hat v^q(\hat N_j)\qquad \forall \hat N_j \in \hat N_i^{p,q}.
  \end{equation}
  We call $\hat R_i^{q,p}$ the anisotropic restriction operator on the reference element $\hat K$ in
the $i$-th direction.
\end{Definition}
\begin{remark}
  The finite elements with anisotropic polynomial degrees are only used in local
  finite element restrictions and interpolations. The global restriction and
  interpolation operators are defined such that they immediately map back
  into finite element spaces with isotropic polynomials. Therefore, we omit the
  explicit definition of the global finite element spaces with anisotropic polynomial degree.
\end{remark}
\begin{Definition}\label{def: aniso space restriction}
  The restriction operator in $i$-th direction $R^p_i:{V}_{ h}^{2p} \to {V}_{ h}^{2p}$ is defined as
  \begin{equation}\label{eq:aniso-space-restriction}
    (R^p_iv_h)(x)\coloneqq \left(\hat R^{2p,p}_i\hat v_h\right)\left(\left(\boldsymbol{T}_K\right)^{-1}(x)\right)\qquad x \in K\,,K \in \mathcal{T}_{h}\,,
  \end{equation}
  for $v_h \in {V}_{ h}^{2p}$ and $\hat v_h =
  \boldsymbol{T}_K^{-1}\circ v_h$.
\end{Definition}
Next, we introduce notation for the a posteriori error representation. Firstly we
define a local anisotropic remainder term and show that it is of higher order,
similar to the remainder terms in the isotropic error estimator~\eqref{eq:3:16}.
\begin{Definition}\label{def: isotrtropic remainder reference}
  Let $\hat v_{2p} \in \hat{Q}_{2p}(\hat{K})$. We define
  $$\hat{\mathbb{E}}^{2p}\hat v_{2p}\coloneqq \hat v_{2p}+\hat R^{2p,p}\hat v_{2p}-\hat R_1^{2p,p}\hat v_{2p}-\hat R_2^{2p,p}\hat v_{2p}$$
  as the local anisotropic remainder on the reference element $\hat{K}$.
\end{Definition}
We now prove, that the anisotropic remainder is in fact of higher order. For the
sake of simplicity and brevity, we restrict ourselves to the case where $p=1$ and $q=2$.
\begin{Lemma}\label{lem: remainder hat K}
  Let $d=2$, $p=1$ and $q=2p=2$. Additionally, let $\hat v_{2p} \in \hat{Q}_{2p}(\hat{K})$.
  For
  $\hat{\mathbb{E}}^{2p}\hat v_{2p}$
  as in Definition~\ref{def: isotrtropic remainder reference}
  holds
  \begin{equation}
    (\hat{\mathbb{E}}^{2p}\hat v_{2p})(\hat{x}_1,\hat{x}_2)=\frac{1}{4}\frac{\partial^4 \hat v_{2p}}{\partial \hat{x}_1^2 \partial \hat{x}_2^2}(\hat{x}_1,\hat{x}_2)\hat{x}_1 \hat{x}_2 (1 - \hat{x}_1) (1 - \hat{x}_2)\,, \qquad (\hat{x}_1, \hat{x}_2) \in \hat{K}.
  \end{equation}
\begin{proof}
  Follows from technical computations with the nodal values and their corresponding functions on $\hat{K}$ and is omitted for the convenience of the reader.
\end{proof}
\end{Lemma}

\begin{Definition} \label{def: global error part}
  Let $v_h \in V_h^{2p}$.   We define
  \begin{equation}
    ({\mathbb{E}}^{2p}v_h)(x)\coloneqq \left(\hat{\mathbb{E}}^{2p}\hat v_h\right)\left(\boldsymbol{T}_K^{-1}(x)\right)\,,\qquad x \in K\,,K \in \mathcal{T}_{h}.
  \end{equation}
  as the isotropic remainder in space.
\end{Definition}
Next we estimate the anisotropic remainder term for arbitrary $K\in \mathcal
T_h$. Here we restrict ourselves to affine mappings $\boldsymbol T_K$ with a
diagonal matrix.
\begin{theorem}\label{thm: bound on Evh}
  Let $p=1$ and $K \in \mathcal{T}_h$ with $K=(\tilde{x}_{1,K},\tilde{x}_{2,K})+[0,h_{1,K}]\times[0,h_{2,K}]$  for some $\tilde{x}_{1,K},\tilde{x}_{2,K}$.
  Furthermore, let $\boldsymbol{T}_K: \hat{K} \to K$ be defined as $\boldsymbol{T}_K \left(\hat{x}_1,\hat{x}_2\right)\coloneqq (x_1,x_2)\coloneqq(\tilde{x}_{1,K}+\hat{x}_1 h_{1,K},\tilde{x}_{2,K}+\hat{x}_2 h_{2,K})$. Then for an arbitrary $v_h \in V_h^{2p}$ holds that
    \begin{equation}
    \left|{\mathbb{E}}^{2p}v_h\right| \leq \frac{1}{16}\sup_{K \in \mathcal{T}_h,(x_1,x_2) \in K}\left|\frac{\partial^4v_{h|K}}{\partial {x}_1^2 \partial {x}_2^2}(x_1,x_2)\right|h_{1,K}^2 h_{2,K}^2.
  \end{equation}

  \begin{proof}
    Let $v_h \in V_h^{2p}$ and $K \in \mathcal{T}_h$. Furthermore, let $\hat v_{2p}=\boldsymbol{T}_K^{-1}\circ \hat v_{2p|K}$.
    Therefore, we get for all  $(\hat{x}_1, \hat{x}_2) \in \hat{K}$ that
    \begin{align*}
      |\hat{\mathbb{E}}^{2p}\hat
      v_{2p}(\hat{x}_1,\hat{x}_2)|&=\frac{1}{4}\left|\frac{\partial^4\hat v_{2p}}{\partial \hat{x}_1^2 \partial \hat{x}_2^2}(\hat{x}_1,\hat{x}_2)\hat{x}_1 \hat{x}_2 (1 - \hat{x}_1) (1 - \hat{x}_2)\right| \\
      &=\frac{1}{4}\left|\frac{\partial^4v_{h|K}}{\partial {x}_1^2 \partial {x}_2^2}(\boldsymbol{T}_K \left(\hat{x}_1,\hat{x}_2\right))h_{1,K}^2 h_{2,K}^2\hat{x}_1 \hat{x}_2 (1 - \hat{x}_1) (1 - \hat{x}_2)\right| \\
      &\leq \frac{1}{4}\left| \frac{\partial^4v_{h|K}}{\partial {x}_1^2 \partial {x}_2^2}(\boldsymbol{T}_K \left(\hat{x}_1,\hat{x}_2\right)) h_{1,K}^2 h_{2,K}^2\right| \frac{1}{4}.
    \end{align*}
    We get that
  \begin{equation*}
    \left|{\mathbb{E}}^{2p}v_h\right| \leq \frac{1}{16}\sup_{K \in \mathcal{T}_h,(x_1,x_2) \in K}\left|\frac{\partial^4v_{h|K}}{\partial {x}_1^2 \partial {x}_2^2}(x_1,x_2)\right|h_{1,K}^2 h_{2,K}^2.
  \end{equation*}
  \end{proof}
\end{theorem}

\begin{remark} \label{remark: also parallelograms} An extension of
  Theorem~\ref{thm: bound on Evh} to parallelograms as well as a higher
  polynomial degree $p>1$ is straight forward.
\end{remark}
\begin{figure}[htb]
  \begin{minipage}[t]{0.49\textwidth}\centering
    {\input{Tikz/FigureQ21.tex} }
    \subcaption{$\Q{2,1}$ finite elements on $K_{2h,1^\ast}^P$. \label{Fig: PatchQ21 Node}}
  \end{minipage}
  \begin{minipage}[t]{0.49\textwidth}\centering
  {\input{Tikz/FigureQ12.tex} }
  \subcaption{$\Q{1,2}$ finite elements on $K_{2h,2^\ast}^P$. \label{Fig: PatchQ12 Node}}
\end{minipage} 
\caption{DoFs of finite elements with anisotropic (a), (b) polynomial degrees on a patch.\label{Fig: Patch QAnIso Nodes}}
\end{figure}
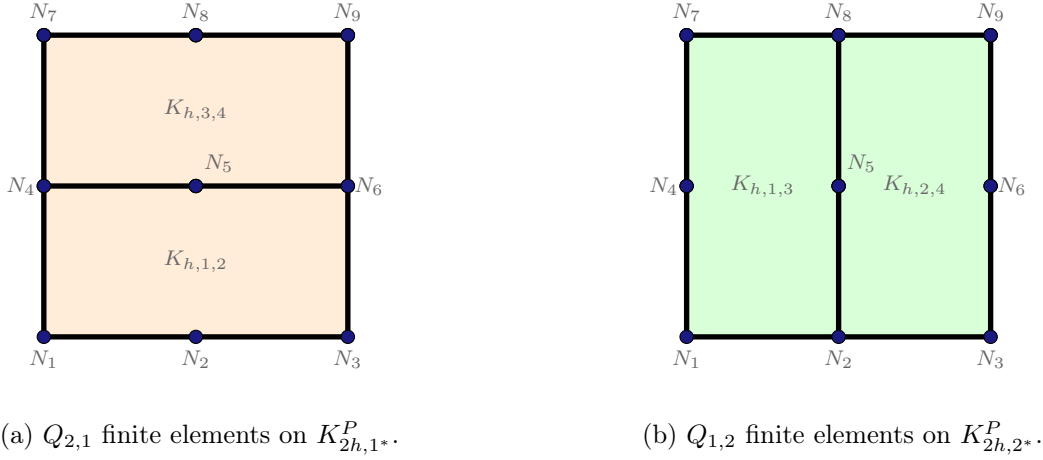

\begin{figure}[htb]
  \begin{minipage}{0.245\textwidth}
    \scalebox{0.6}{\input{Tikz/PatchQ11HighInter.tex} }
    \subcaption{$v_h$.}
  \end{minipage} \hfill%
  \begin{minipage}{0.245\textwidth}
    \scalebox{0.6}{\input{Tikz/PatchQ21HighInter.tex} }
    \subcaption{$\mathrm{I}_{2h,1}^{\left(2p\right)}v_h$.}
  \end{minipage} \hfill%
  \begin{minipage}{0.245\textwidth}
    \scalebox{0.6}{\input{Tikz/PatchQ12HighInter.tex}}
    \subcaption{$\mathrm{I}_{2h,2}^{\left(2p\right)}v_h$.}
  \end{minipage}\hfill%
  \begin{minipage}{0.245\textwidth}
    \scalebox{0.6}{\input{Tikz/PatchQ22HighInter.tex}}
    \subcaption{$\mathrm{I}_{2h}^{\left(2p\right)}v_h$.}
  \end{minipage}%
  \caption{The action of different interpolation operators on a function $v_h$ on the path $K_{2h}$ and $p=1$.\label{Fig: Action AnIso Int}}
\end{figure}
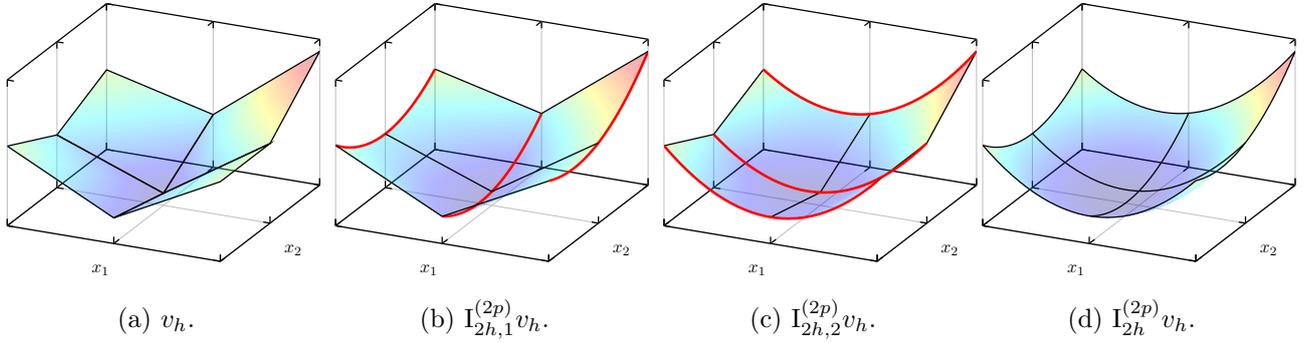
\begin{Definition} \label{def: AnIsoHighInt}
  For $K_{2h}^P\in \mathcal{T}_{2h}^P$, $K_{2h}\in \mathcal{T}_{2h}$, we define
  $\mathrm{I}_{2h,i}^{\left(2p\right)}\colon V_h^{p}(K_{2h}^P) \to V_h^{2p}(K_{2h})$, the patch-wise higher order interpolation into the spatial direction $i$, as
  \begin{equation}
    \mathrm{I}_{2h,i}^{\left(2p\right)}\coloneqq R^{p}_i\circ\mathrm{I}_{2h}^{\left(2p\right)}
  \end{equation}
\end{Definition}

With abuse of notation, the operator $R_i^p$
(cf.~\eqref{eq:aniso-space-restriction}), which was originally defined for
elements $K \in \mathcal{T}_h$, is here extended to act on elements $K_{2h} \in
\mathcal{T}_{2h}$. Further, the global anisotropic restriction operator is denoted by
the same name. The action of directional interpolation operators, where $p=1$, is visualized in
Figure~\ref{Fig: Action AnIso Int} and the DoFs are shown in Figure~\ref{Fig: Patch QAnIso Nodes}.

\subsection{Anisotropic Goal Oriented Error Estimation}
In this subsection, we extend the classical error representation from
Theorem~\ref{Thm:3.3} to the anisotropic setting, following the approach
introduced by Richter in~\cite[Theorem~2]{richter_posteriori_2010}. The key
idea is the anisotropic splitting of the discretization error into distinct
directional contributions, providing the theoretical foundation for anisotropic
goal-oriented error estimation. While our framework closely
follows the one developed in~\cite{richter_posteriori_2010}, our contribution is a different
approximation of the directional error contributions, which we precisely state
in the following theorem.

  \begin{theorem}\label{thm:anisotropic-error-split}
    Let $\eta_{h}$ be defined as in \eqref{eq: implemented error estimator}.
    Then, for the splitting into the single directions it holds
    \begin{equation}
        \eta_{h}\coloneqq\eta_{h,1}+ \eta_{h,2}+\eta_{h,\mathbb{E}},
    \end{equation}
  where
    \begin{equation}
      \begin{aligned}\label{eq:aniso-error-est-h-aniso}
      \eta_{h,i}&\coloneqq\frac{1}{2}\rho_{\tau}(u_{\tau h})(z_{\tau h}- R^{p}_iz_{\tau h})
      + \frac{1}{2}
      \rho_{\tau}^{\ast}(u_{\tau h},z_{\tau h})
      (\mathrm{I}_{2h,i}^{\left(2p\right)}u_{\tau h}-u_{\tau h})
      \\
      & \qquad
      + \frac{1}{2} S(u_{\tau h})( R^{p}_iz_{\tau h})
      + \frac{1}{2} S^{\prime}(u_{\tau h})(\mathrm{I}_{2h,i}^{\left(2p\right)}u_{\tau h}-u_{\tau h},\mathrm{R}_{h}^{p}z_{\tau h})\,,
    \end{aligned}
  \end{equation}
for $i \in \{1,2\}$ and
    \begin{equation}
  \begin{aligned}
    \eta_{h,\mathbb{E}}&\coloneqq\frac{1}{2}\rho_{\tau}(u_{\tau h})(\mathbb{E}^{2p}z_{\tau h})
    + \frac{1}{2}
    \rho_{\tau}^{\ast}(u_{\tau h},\mathrm{R}_{h}^{p}z_{\tau h})
    (\mathbb{E}^{2p}\mathrm{I}_{2h}^{\left(2p\right)}u_{\tau h})
    \\
    & \qquad
    + \frac{1}{2} S(u_{\tau h})(\mathbb{E}^{2p}z_{\tau h})
    + \frac{1}{2} S^{\prime}(u_{\tau h})(\mathbb{E}^{2p}\mathrm{I}_{2h}^{\left(2p\right)}u_{\tau h},\mathrm{R}_{h}^{p}z_{\tau h})\,.
  \end{aligned}
\end{equation}
\begin{proof}
  For $v_h \in \mathcal{V}_{h}^{2p}$, we obtain
  from the definition of the restrictions~\eqref{eq:iso-space-restriction} and~\eqref{eq:aniso-space-restriction}, that
  \begin{equation}\label{eqn: help anisthm}
    (I + R^p_h - R^p_1 -R^p_2)v_{h}-\mathbb{E}^{\left(2p\right)}v_{h} = 0\,,
  \end{equation}
  where $I$ is the identity mapping.
  Adding $(I - R^p_h)v_{h}$ on both sides of~\eqref{eqn: help anisthm}, we get
  \begin{equation}\label{eqn: difference approximation}
   v_{h}-R^p_hv_{ h}=\left(I - R^p_1\right)v_{ h} +\left(I- R^p_2\right)v_{ h}-\mathbb{E}^{\left(2p\right)}v_{ h} \,.
\end{equation}
Since by definition $\mathrm{I}_{2h}^{\left(2p\right)}R^p_h=I$, it directly
follows that for  $v_h\in \mathcal{V}_{h}^{2p}$ we have
\begin{equation} \label{eqn: difference approx inter}
   \mathrm{I}_{2h}^{\left(2p\right)}v_h-\underbrace{R^p_h\mathrm{I}_{2h}^{\left(2p\right)}v_h}_{=v_h}=\left(I - R^p_1\right)\mathrm{I}_{2h}^{\left(2p\right)}v_h +\left(I- R^p_2\right)\mathrm{I}_{2h}^{\left(2p\right)}v_h-\mathbb{E}^{\left(2p\right)}\mathrm{I}_{2h}^{\left(2p\right)}v_h \,.
\end{equation}
Furthermore, from~\eqref{eqn: help anisthm} we deduce that
\begin{equation}\label{eqn: sum approximation}
  (I + R^p_h)v_{h} =  (R^p_1 +R^p_2)v_{h}+\mathbb{E}^{\left(2p\right)}v_{h}\,.
\end{equation}
The discrete isotropic error estimator $\eta_{h}$, given in \eqref{eq:
  implemented error estimator}, then reads as
\begin{equation*}
  \begin{aligned}
    \eta_{h}&=\frac{1}{2}\rho_{\tau}(u_{\tau h})(z_{\tau h}-\mathrm{R}_{h}^{p}z_{\tau h})
    + \frac{1}{2}
    \rho_{\tau}^{\ast}(u_{\tau h},\mathrm{R}_{h}^{p}z_{\tau h})
    (\mathrm{I}_{2h}^{\left(2p\right)}u_{\tau h}-u_{\tau h})
    \\
    & \qquad
    + \frac{1}{2} S(u_{\tau h})(z_{\tau h}+\mathrm{R}_{h}^{p}z_{\tau h})
    + \frac{1}{2} S^{\prime}(u_{\tau h})(\mathrm{I}_{2h}^{\left(2p\right)}u_{\tau h}-u_{\tau h},\mathrm{R}_{h}^{p}z_{\tau h})\,.
  \end{aligned}
\end{equation*}
We obtain the result by substiting  $z_{\tau h}-\mathrm{R}_{h}^{p}z_{\tau h}$,
$\mathrm{I}_{2h}^{\left(2p\right)}u_{\tau h}-u_{\tau h}$ and $z_{\tau
  h}+\mathrm{R}_{h}^{p}z_{\tau h}$ by the equations~\eqref{eqn: difference
  approximation} to~\eqref{eqn: sum approximation}, respectively.
\end{proof}
  \end{theorem}

\begin{remark}\label{rem:localization}
  Here $\eta_{h,1}$ and $\eta_{h,2}$ are the directional error
  estimators. For localization, we employ the classical approach as in~\cite{BeRa01,BaRa03}, analogous to the isotropic error estimator given in~\eqref{eq: implemented error estimator}.
  A partition of unity technique as in \cite{RiWi15_dwr,endtmayer_goal-oriented_2024,ThiWi24,EnLaWi20,beuchler2023multigoal,beuchler2024mathematical} is possible as well.
 \end{remark}
\begin{remark}\label{rem:negelecting}
  Since, the term
  $\mathbb{E}^{2p} v_h$ is of higher order
  (cf.~Theorem~\ref{thm: bound on Evh}), we omit it in our numerical
  examples and do not restrict ourselves to parallelograms in the
  experiments (cf.~Remark~\ref{remark: also parallelograms}).
\end{remark}

\section{Algorithm for Goal-Oriented Anisotropic Mesh Adaptivity}\label{sec:algorithm}
The element wise contributions of the anisotropic error indicators in space and time are given by
\begin{subequations}\label{eq:localized-aniso-error-ind}
\begin{align}\label{eq:spatial_localized_anisotropic_indicators}
\eta_{\tau}&=\sum_{n=1}^{N}\sum_{K\in\mathcal{T}_h^n}\eta_{\tau}^{K,n}\,,\\
\label{eq:temporal_localized_anisotropic_indicators}
\eta_{h}&=\sum_{i=1}^d\sum_{n=1}^{N} \sum_{K\in\mathcal{T}_h^n} \eta_{h,i}^{K,n}\,,
\end{align}
\end{subequations}
where we neglect the spatial higher order remainder term of $\eta_{h}$ in Theorem~\ref{thm:anisotropic-error-split} according to Rem.~\ref{rem:negelecting} and $\eta_{\tau}$ is defined in~\eqref{eq: implemented error estimator-tau}, respectively.
These indicators steer the anisotropic mesh
adaptation strategy presented in Algorithm~\ref{alg:AnisotropicAMR}. There, we
use the spatial directional indicators and the temporal indicators to mark
spatial elements in each direction and time intervals. This marking results in
anisotropic and isotropic refinements.
  \begin{algorithm}[htb]
    \caption{\label{alg:AnisotropicAMR}Anisotropic Mesh Adaptation: mark \& refine}
    \begin{algorithmic}[1]
      \REQUIRE{$\displaystyle\eta^{K,n}_{h,i}$,
        $\eta^{K,n}_{\tau}$ (cf.~\eqref{eq:localized-aniso-error-ind}) for all $K\in
        \mathcal T_h$, for all $I_n\in \mathcal T_{\tau}$, $\theta_h$, $\theta_\tau$}
      \STATE{\textbf{Calculate} for all $\displaystyle K\in \mathcal T_{h},\:i=1,\dots,d$:
        $\displaystyle\eta^{K}_{h,i}=\displaystyle\sum_{I_n\in \mathcal T_{\tau}}\eta^{K,n}_{h,i}$.}\label{alg:calc-eta-khi}
      \STATE{\textbf{Calculate} for all $\displaystyle I_n\in \mathcal T_{\tau}$:
        $\displaystyle\eta^{n}_{\tau}=\sum_{K\in \mathcal T_{h}}\eta^{K,n}_{\tau}$.}\label{alg:calc-eta-nt}
      \FOR{$i=1,\dots,d$}
      \STATE{\textbf{Mark} $\displaystyle\theta_h \lvert\mathcal T_{h}\rvert$
        cells $\displaystyle K\in \mathcal
        T_h$ with the largest refinement criteria $\displaystyle \eta^{K}_{h,i}$ for \emph{refinement in direction} $i$.}\label{alg:mark-aniso}
      \ENDFOR{}
      \STATE{\textbf{Mark} $\displaystyle\theta_\tau \lvert\mathcal T_{\tau}\rvert$ subintervals $\displaystyle I_n\in \mathcal
        T_\tau$ with the largest refinement criteria $\displaystyle \eta^{n}_{\tau}$ for refinement.}\label{alg:mark-time}
      \STATE{\textbf{Refine} in space according to marking in line~\ref{alg:mark-aniso}}
      \STATE{\textbf{Refine} in time according to marking in line~\ref{alg:mark-time}}
    \end{algorithmic}
  \end{algorithm}
  \begin{remark}[Anisotropic Mesh Adaptation]
  ~\\[-2\topsep]
  \begin{itemize}\itemsep1pt \parskip0pt \parsep0pt
  \item As the spatial mesh is fixed for all time intervals, the local error
    indicators \(\eta^{K,n}_{h,i}\), computed for each cell \(K\) and direction
    \(i\), are summed over time intervals to yield the total directional error
    \(\eta^K_{h,i}\) (Algorithm~\ref{alg:AnisotropicAMR},
    line~\ref{alg:calc-eta-khi}).
    The temporal error \(\eta^{n}_{\tau}\) results from summing local
    contributions across all spatial elements
    (Algorithm~\ref{alg:AnisotropicAMR}, line~\ref{alg:calc-eta-nt}).
  \item Marking for refinement is performed separately in each spatial direction
    (Algorithm~\ref{alg:AnisotropicAMR}, line~\ref{alg:mark-aniso}). A spatial
    element \(K\) is marked for refinement in direction \(i\) if its directional
    indicator \(\eta^K_{h,i}\) is among the top \(\theta_h |\mathcal T_{h}|\)
    values. Evaluating directional errors individually naturally yields anisotropic
    refinement when one direction dominates, and isotropic refinement when
    errors are balanced.
  \item Adaptive refinement on quadrilateral or hexahedral meshes leads to
    interdependent hanging nodes. To resolve these mutual dependencies, we refine the
    coarser element at the end of the hanging-node chain in the appropriate
    direction (i.\,e.\ anisotropically).
  \end{itemize}
  \end{remark}
  To compute the error indicators required in
  Algorithm~\ref{alg:AnisotropicAMR}, we first solve the primal problem forward
  in time, followed by the adjoint problem backward in time. Using both primal
  and adjoint solutions, we calculate the local error estimators~\eqref{eq:localized-aniso-error-ind}, perform
  marking and refinement according to Algorithm~\ref{alg:AnisotropicAMR}, and
  repeat the space-time solution procedure iteratively until the desired
  accuracy is reached. For algorithmic details and implementation aspects,
  especially regarding space-time finite elements, we refer
  to~\cite{KoeBruBau2019,BruchhaeuserB22,BruchhaeuserBR12,BruchhaeuserSV08,RoThiKoeWi23,ThiWi24,endtmayer_goal-oriented_2024,dorfler_space-time_2016,dorfler_petrovgalerkin_2023}.

  \section{Numerical Examples}\label{sec:8:numerical_examples}
  To validate the efficacy of our anisotropic adaptive mesh refinement approach applied to time-dependent CDR equations, we present three numerical examples.
  \begin{enumerate}
  \item \emph{Interior Layer Problem:}
  For this benchmark case an exact solution is available, which is characterized by a sharp interior layer.
  \item \emph{Stationary Hemker Problem:} The Hemker problem serves as a more challenging test case. It can be interpreted as a model of a convection-dominated
    heat transfer from a hot column. The solution to the Hemker problem is
    characterized by a boundary layer at the cylinder's surface and two interior
    layers located downstream of the cylinder.
  \item \emph{Nonstationary Hemker Problem with quadratic obstacle:}
    We modify the classical Hemker problem to be nonstationary and replace the
    circular obstacle by a quadratic one.
  \end{enumerate}
  In these examples, we employ isotropic and anisotropic refinement strategies.
  Depending on the refinement type, a different marking strategy is employed. In
  the case of isotropic refinement, dynamic meshes are used, i.\,e.\ the mesh
  can be different between timesteps. Thereby, each time slab has its own
  mesh. Conversely, the anisotropic mesh refinement strategy only allows fixed
  meshes in time at the current stage of development.

  To measure the accuracy of the error estimator we will study the
  \emph{effectivity index}
  \begin{equation}
    \Ieff^{\textup{a}} = \left\vert\frac{\eta_{h,\,x}+\eta_{h,\,y}+\eta_{\tau}}{J(u) - J(u_{\tau h})}\right\vert,\quad
    \Ieff = \left\vert\frac{\eta_{h}+\eta_{\tau}}{J(u) - J(u_{\tau h})}\right\vert\,,
  \end{equation}
  for the anisotropic and isotropic case, respectively. For the anisotropic
  case, we further denote
  $\eta_{h}^{\textup{a}}=\eta_{h,\,x}+\eta_{h,\,y}$ and
  $\eta_{\tau h}^{\textup{a}}=\eta_{\tau}+\eta_{h}^{\textup{a}}$, where $\eta_{h,\,x} \coloneqq \eta_{h,1}$ and  $\eta_{h,\,y} \coloneqq \eta_{h,2}$. For the
  isotropic case we denote $\eta_{\tau h}=\eta_{\tau}+\eta_{h}$.
  Moreover, as an indicator for the anisotropy we consider the maximum
  \emph{aspect ratio}, given by
  \begin{equation}
    \operatorname{ar}_{\max}\coloneqq \max_{K \in \mathcal T_h} \max_{\boldsymbol x \in
      Q_K}\frac{\lambda_{\max}}{\lambda_{\min}}\,,
  \end{equation}
  where $Q_K$ is the set of quadrature points on $N_{\text{time}}$ and $\lambda_{\min},
  \lambda_{\max}$ are the minimal and maximal eigenvalues of \(\nabla
  \boldsymbol T_K(\boldsymbol x)\), respectively.
  The SUPG stabilization paramater $\delta_{K}(\cdot)$  is defined as
  \[
    \delta_K = \delta_0\,h_K,
  \]
  where $h_K=\sqrt[d]{|K|}$ is the cell diameter of the mesh cell. Throughout this work we
  choose $\delta_0=0.1$. The total number of DoFs  is denoted by $N_{\text{tot}}$, whereas the spatial DoFs and temporal DoFs are denoted by $N_{\text{space}}$ and $N_{\text{time}}$, respectively.
   The implementation is based on the \texttt{deal.II}
  finite element library~\cite{africa_dealii_2024}. The tests are run on a
  single node with 2 Intel Xeon Platinum 8360Y CPUs and \SI{1024}{\giga\byte}
  RAM of the HPC cluster HSUper at HSU.\@

  \subsection{Interior Layer Problem}
  \label{sec:8.1:step-layer}
  This well-known benchmark with a sharp interior layer of thickness
  $\mathcal{O}(\sqrt{\varepsilon}|\log \varepsilon|)$ has an exact solution
  \begin{equation}
    \label{eq:step-layer}
    u(\mathbf{x},\,t) = \frac{\operatorname{e}^{3(t-1)}}{2}\left( 1 - \tanh \frac{2x-y-\frac{1}{2}}{\sqrt{5\varepsilon}}\right)\,.
  \end{equation}
  The problem is defined on
  $\mathcal{Q}\coloneqq\Omega\times I=(0,\,1)^2\times (0,\,1]$ with inhomogeneous boundary
  conditions given by~\eqref{eq:step-layer} on the whole spatial boundary
  $\Gamma_D=\partial \Omega$. The convection is set to
  $\mathbf{b}=\frac{1}{\sqrt{5}}(1,\,2)^{\top}$, the diffusion to
  $\varepsilon=\num{1.e-6}$ and the reaction coefficient is given by $\alpha=1$.
  As sketched in Figure~\ref{fig:interior-layer-geo}, we use structured and
  unstructured meshes with isotropic and anisotropic refinement. In all
  configurations, we refine and coarsen a fixed fraction of cells in both space
  and time. In space, the refinement fraction is set to
  \(\theta_{\text{space}}^{\text{ref}} = \frac{1}{5}\), and the coarsening
  fraction is \(\theta_{\text{space}}^{\text{co}} = \frac{1}{100}\). Similarly,
  in time, the refinement fraction is
  \(\theta_{\text{time}}^{\text{ref}} = \frac{2}{3}\), while no coarsening is
  applied, i.e., \(\theta_{\text{time}}^{\text{co}} = 0\).
  Finally, the goal functional is chosen to control the global $L^2(L^2)$-error in space and time, given by
\begin{equation}
\label{eq:42:JL2L2}
J(u)= \frac{1}{\|e\|_{\mathcal{Q}}}\displaystyle\int_I(u,e)\mathrm{d}t\,,
\quad \mathrm{with} \;\; \|\cdot\|_{\mathcal{Q}}
= \left(\int_I(\cdot,\cdot)\;\mathrm{d}t\right)^{\frac{1}{2}},\; e:=u-u_{\tau h}\,.
\end{equation}
  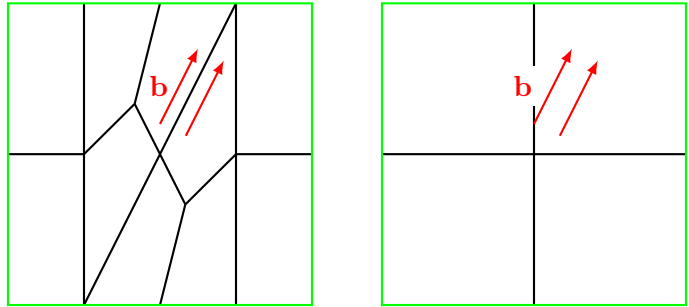
\begin{figure}[!htb]
  \centering
    \begin{tikzpicture}[scale=4,thick]
      \draw[black] (0.25, 0) -- ++ (0,1);
      \draw[black] (0, 0.5) -- ++ (0.25,0);
      \draw[black] (0.75, 0) -- ++ (0,1);
      \draw[black] (0.75, 0.5) -- ++ (.25,0);
      \coordinate (A) at (0.4166666,0.6666666);%
      \coordinate (B) at (0.5833333,0.3333333);%
      \draw[black] (A) -- (B);
      \draw[black] (0.25, 0) -- (0.75, 1);
      \draw[red,-latex] (0.5, 0.6) --node[left,red]{$\mathbf{b}$} ++ (0.125,0.25);
      \draw[red,-latex] (0.585, 0.561) -- ++ (0.125,0.25);
      \draw[black] (A) -- (0.25, 0.5);
      \draw[black] (A) -- (0.5, 1);
      \draw[black] (B) -- (0.75, 0.5);
      \draw[black] (B) -- (0.5, 0);
      \draw[green] (0, 0) rectangle (1,1);
    \end{tikzpicture}\hspace*{2em}
    \begin{tikzpicture}[scale=4,thick]
      \draw[black] (0.5, 0) -- ++(0,1);
      \draw[black] (0, 0.5) -- ++(1,0);
      \draw[red,-latex] (0.5, 0.6) --node[left=3.75pt,red,fill=white]{$\mathbf{b}$} ++ (0.125,0.25);
      \draw[red,-latex] (0.585, 0.561) -- ++ (0.125,0.25);
      \draw[green] (0, 0) rectangle (1,1);
    \end{tikzpicture}
    \caption{Geometry and coarse, unstructured (left) and structured (right) spatial mesh of the domain $\Omega$ for Example~\ref{sec:8.1:step-layer}.}\label{fig:interior-layer-geo}
  \end{figure}

  \begin{figure}[ht]
    \centering
    \begin{tikzpicture}
    \begin{axis}[%
      width=10cm,%
      height=7.5cm,%
      xlabel={\# number of DoFs},%
      ylabel={$\|e\|_{\mathcal{Q}}$},%
      xmode=log,%
      ymode=log,%
      ]%
      \addplot+%
      table [%
      x index=1,%
      y index=2,%
      col sep=space ] {data/unstructured_aniso/convergence.txt};%
      \label{plot:line1}%
      \addplot+%
      table [%
      x index=1,%
      y index=2,%
      col sep=space ] {data/unstructured_iso/convergence.txt};%
      \label{plot:line2}%
      \addplot+%
      table [%
      x index=1,%
      y index=2,%
      col sep=space ] {data/structured_aniso/convergence.txt};%
      \label{plot:line5}%
      \addplot+%
      table [%
      x index=1,%
      y index=2,%
      col sep=space ] {data/structured_iso/convergence.txt};%
      \label{plot:line4}%
      \addplot+%
      table [%
      x index=1,%
      y index=2,%
      col sep=space ] {data/unstructured_iso_uniform/convergence.txt};%
      \label{plot:line3}%
      \coordinate (legend) at (axis description cs:1.025,0.5);%
    \end{axis}%
    \matrix (legend-m) [%
    align=left,%
    draw=none,%
    matrix of nodes,%
    anchor=west,row sep=1pt,%
    column 1/.style={anchor=base east},%
    column 2/.style={anchor=base west,column sep=0.5cm},%
    column 3/.style={anchor=base west},%
    ] at (legend) {%
      \vphantom{{\bfseries Mesh}}\phantom{\ref{plot:line1}}& {\bfseries Mesh} &
      {\bfseries Refinement} \\ 
      \ref{plot:line1} & Unstructured & Anisotropic \\%
      \ref{plot:line5} & Structured & Anisotropic \\%
      \ref{plot:line2} & Unstructured & Isotropic \\%
      \ref{plot:line4}\vphantom{Isotropic} & Structured & Isotropic \\%
      \ref{plot:line3} & Unstructured & Uniform \\%
    };%
    \draw[thick] ([xshift=-5pt]legend-m-1-1.north west) --
    ([xshift=5pt]legend-m-1-3.north east); 
    \draw ([xshift=-5.5pt]legend-m-1-1.south west) --
    ([xshift=5.5pt]legend-m-1-3.south east); 
    \draw[thick] ([xshift=-5pt]legend-m-6-1.south west) --
    ([xshift=5pt]legend-m-6-3.south east-|legend-m-1-3.north east); 
  \end{tikzpicture}
  \caption{\label{fig:step-layer-error}A comparison of the different refinement
    strategies and meshes in terms of their efficiency measured by the error
    over the number of space-time DoFs.}
  \end{figure}
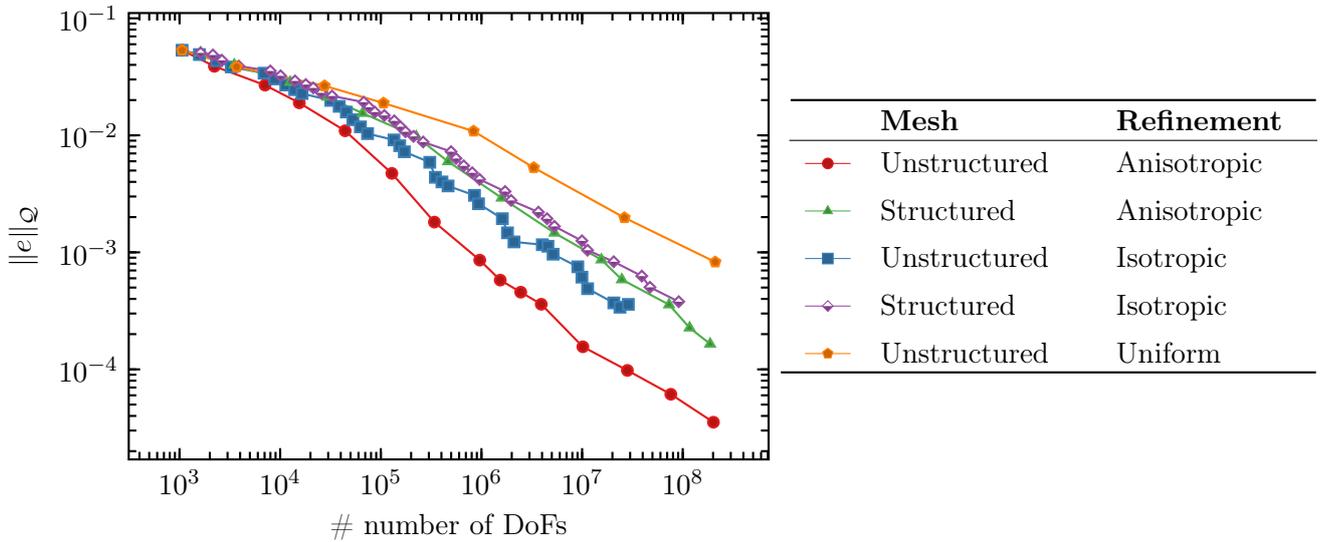
  \begin{figure}[ht]
    \centering
    \begin{tabular}{ccc}
    \includegraphics[width=0.31\textwidth]{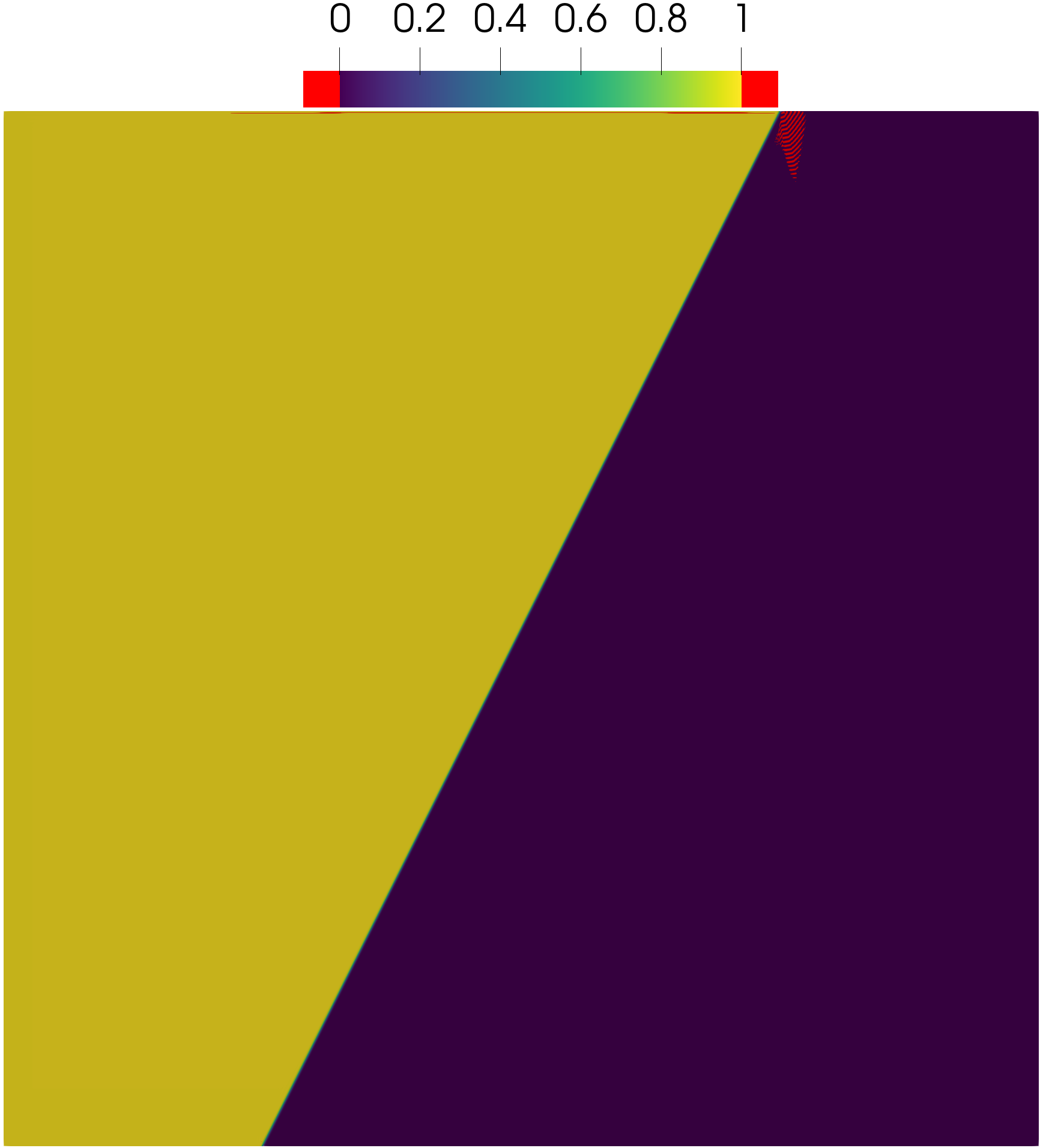} &
    \includegraphics[width=0.31\textwidth]{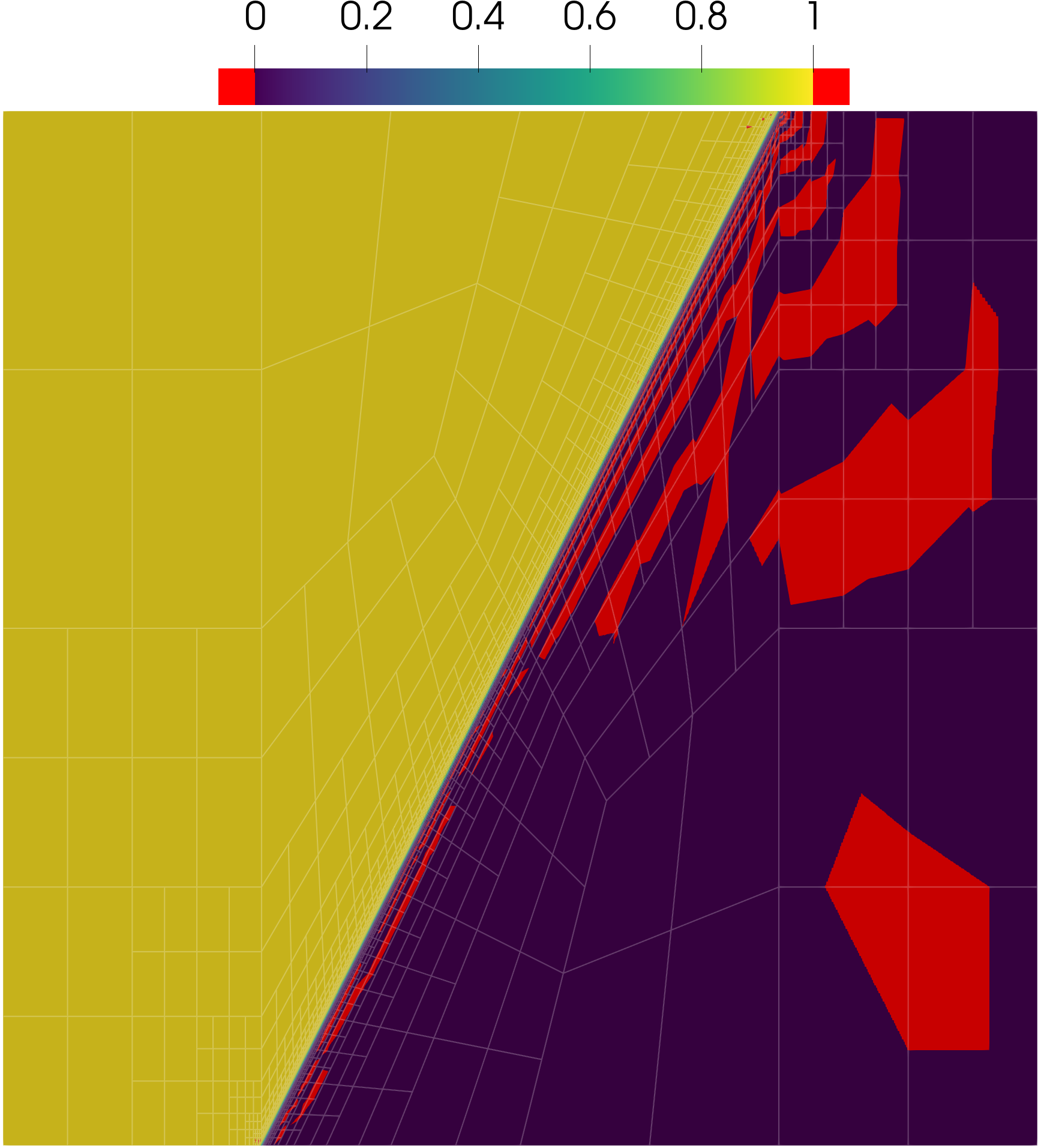} &
    \includegraphics[width=0.31\textwidth]{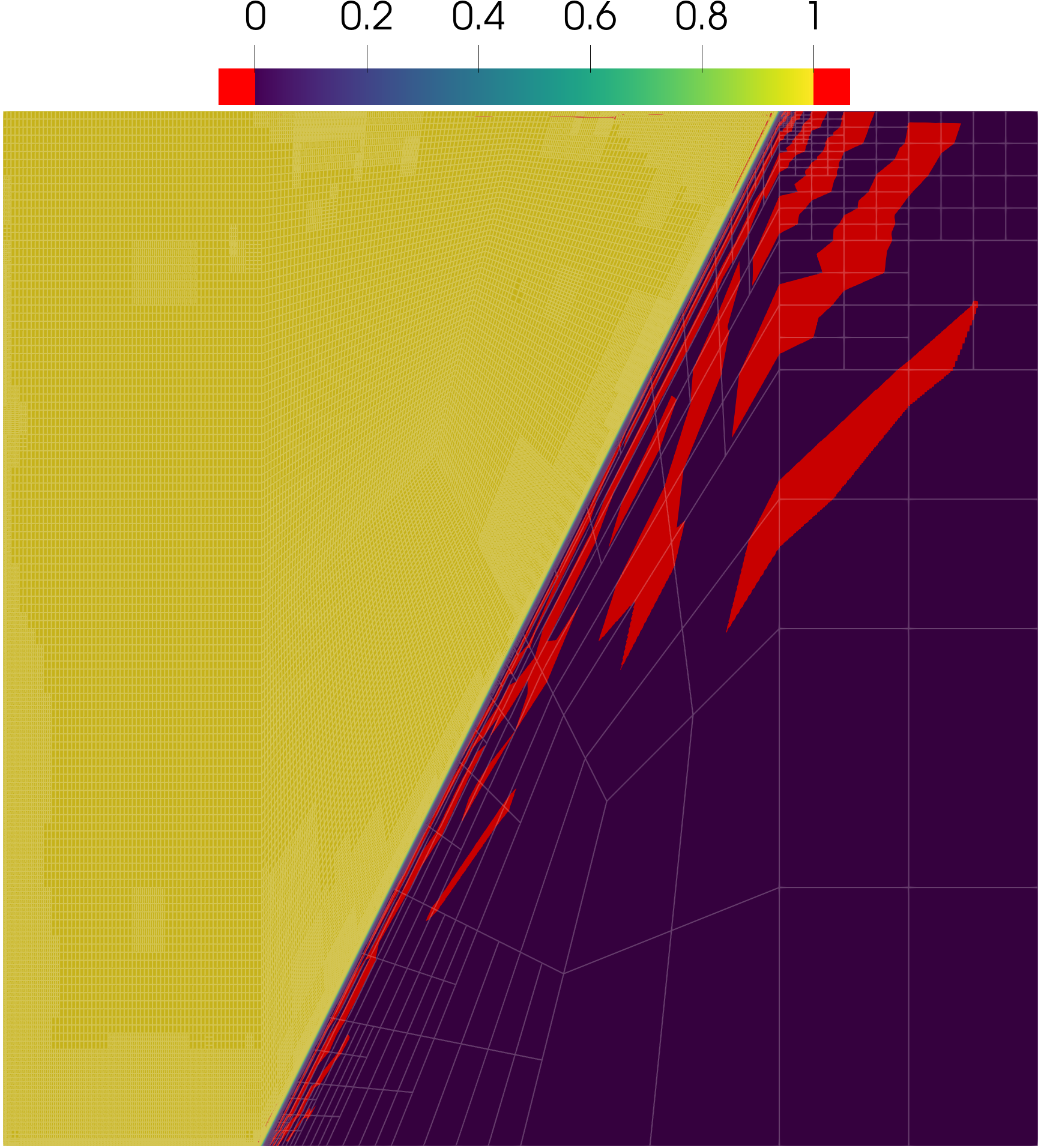}\\
    \includegraphics[width=0.31\textwidth]{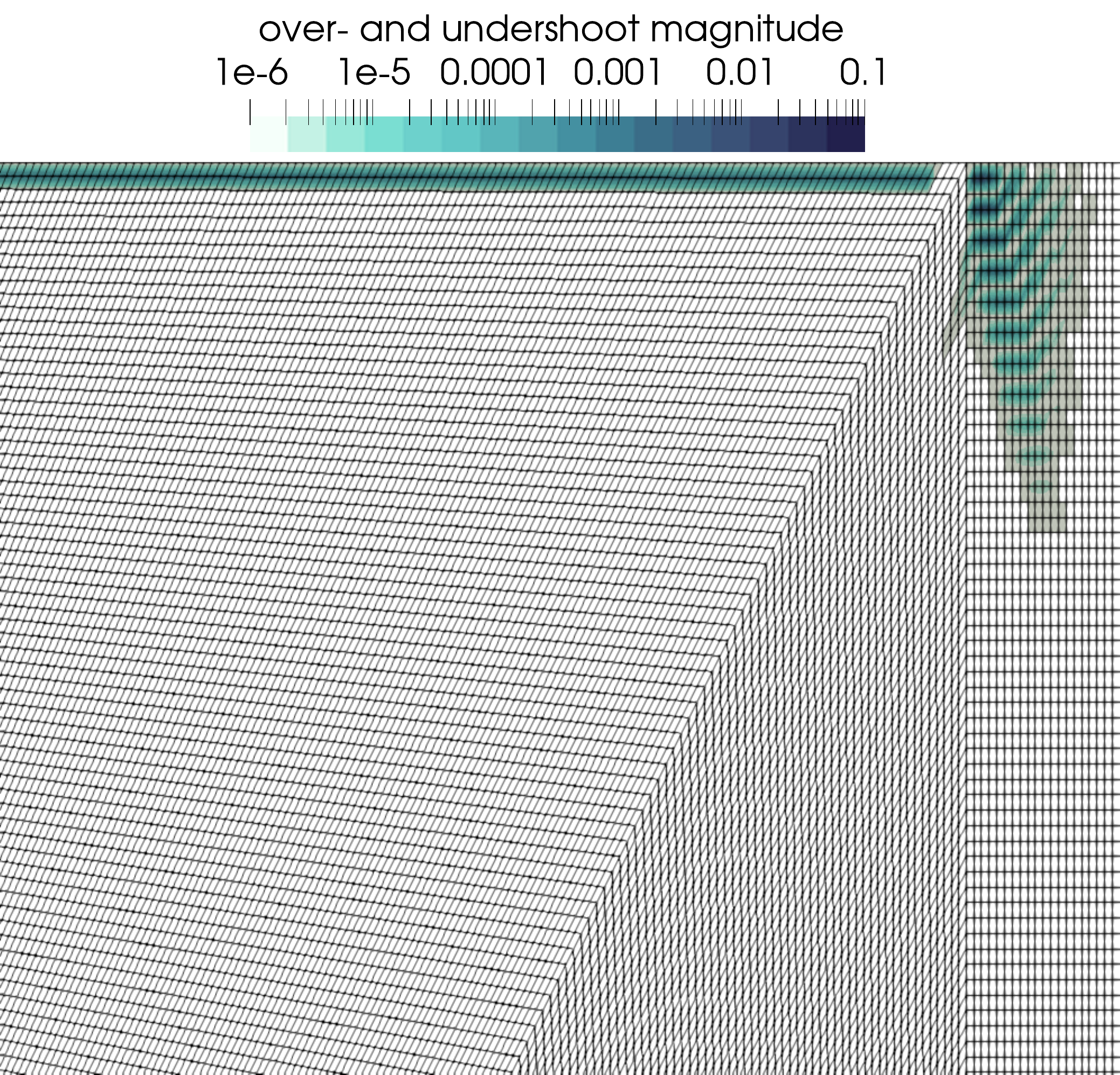} &
    \includegraphics[width=0.31\textwidth]{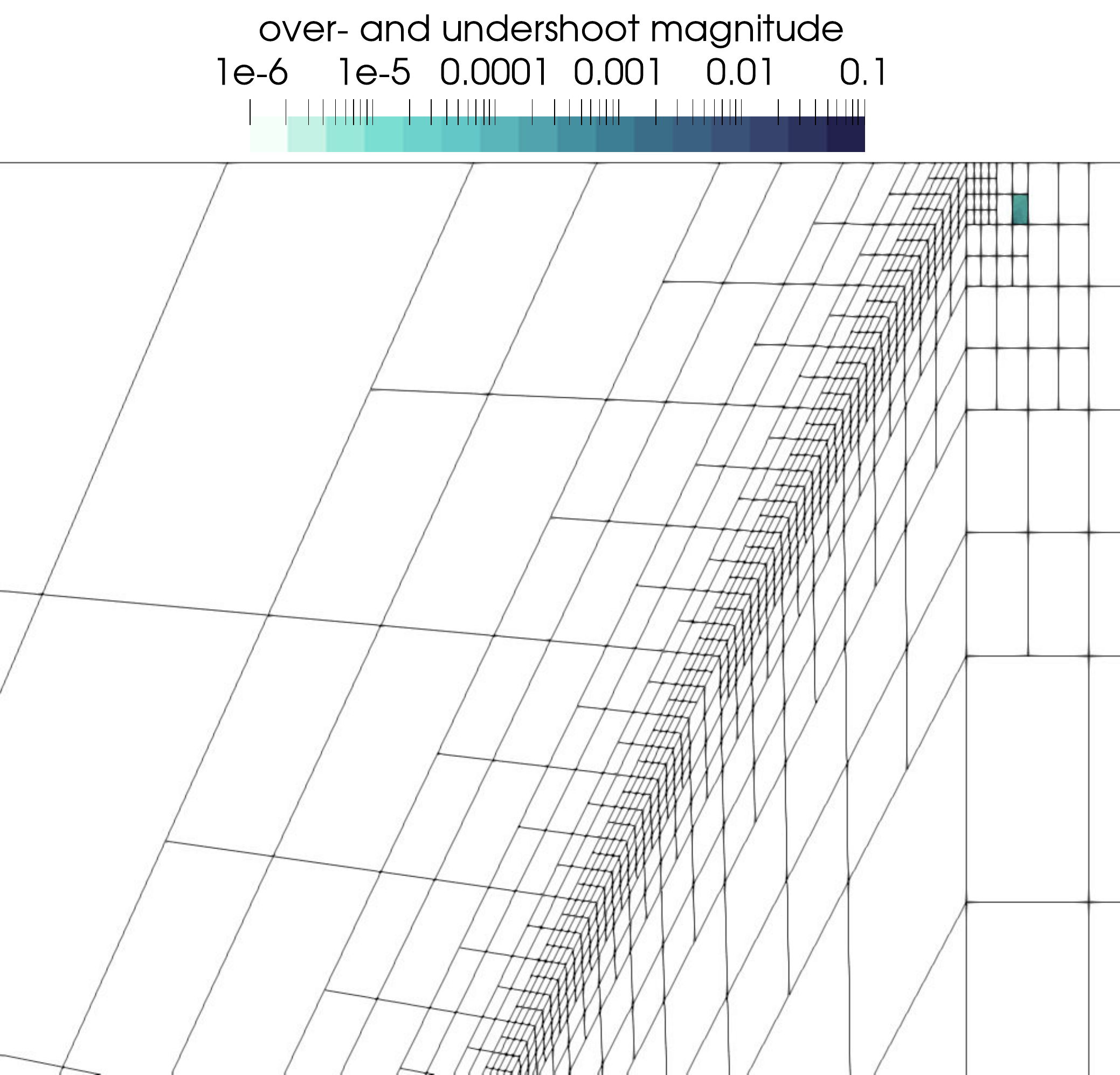} &
    \includegraphics[width=0.31\textwidth]{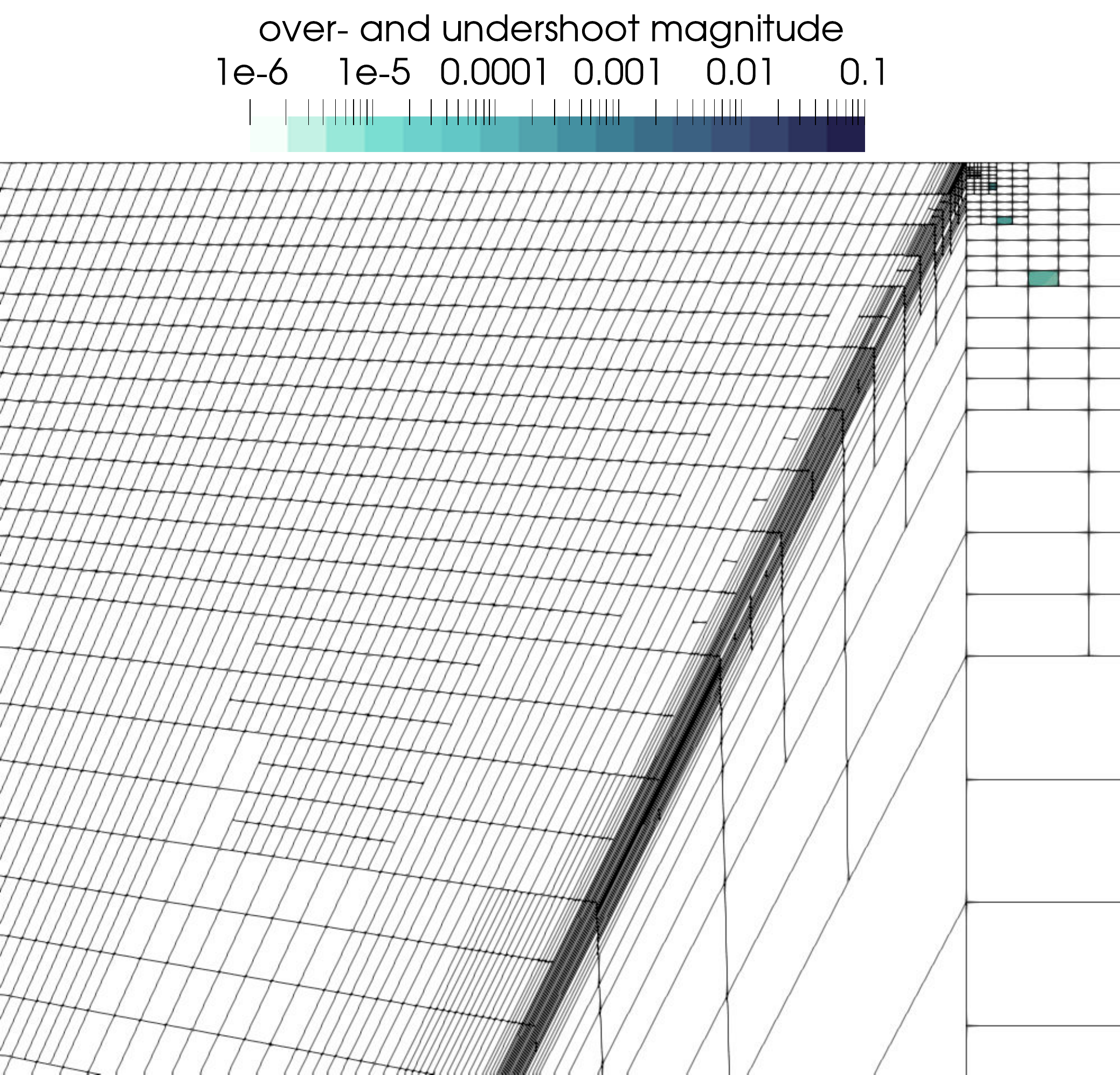} \\
    \end{tabular}
    \caption{All solutions obtained using three different approaches: uniform
      refinement, isotropic adaptive refinement, and anisotropic adaptive refinement.
      Over- and undershoots are highlighted in red, independent of their
      magnitude. The second row illustrates the
      magnitude of the over- and undershoots in the subdomain
      $[0.63,\,0.88] \times [0.77,\,1.0]$. Only regions where the over- and
      undershoot is larger than $\num{1e-6}$ are plotted.}\label{fig:solution_grid}
  \end{figure}
  Figure~\ref{fig:step-layer-error} compares the convergence behavior of the
  error norm \(\|e\|_{\mathcal{Q}}\) versus the number of DoFs for various mesh refinement strategies. The anisotropic
  and isotropic refinements are clearly superior to the uniform refinement on an unstructured mesh. In addition, the anisotropic refinement on unstructured and structured meshes exhibits a faster error reduction per DoF than their isotropic counterparts. Overall, we confirm that anisotropic adaptive refinement achieves superior accuracy with fewer DoFs, which is crucial for efficiently resolving convection-dominated phenomena.

  In Figure~\ref{fig:solution_grid}, we present the solution profiles and spatial
  meshes at final time $T=1$ regarding the different refinement strategies on an unstructured mesh.
  With regard to over- and undershoots, which is an important measure for spurious oscillations for convection-dominated problems, we obtain the following results.

  \paragraph*{Over- and Undershoots} We measure the over- and undershoots in the
  last DWR loop at
  final time $T=1$, where the solution lies within the range \([0,\,1]\). Uniform refinement yields an undershoot of
  \(0.0153\) and an overshoot of \(0.0044\) at a maximum of $\num{656897}$
  spatial DoFs. Isotropic adaptive refinement
  reduces them to \(0.0063\) and \(0.0074\), respectively, at a maximum of $\num{43143}$
  spatial DoFs. The most favorable
  results are obtained with anisotropic adaptive refinement, which produces
  an undershoot of \(7.329\times 10^{-4}\) and an overshoot of
  \(2.9397\times 10^{-4}\) at a maximum of $\num{43143}$
  spatial DoFs.
  Comparing the spatial meshes in Figure~\ref{fig:solution_grid}, we observe a
  significantly higher concentration of elements at the interior layer for the
  anisotropic adaptive refinement compared to the other strategies.

  In Tables~\ref{tab:sl-us-aniso}--\ref{tab:sl-s-iso} we compare the errors and
  effectivity indices for the different refinement strategies.
  In Table~\ref{tab:sl-us-aniso} (unstructured
  anisotropic), the error norm \(\|e\|_{(0,T)\times\Omega}\) decreases to
  \(3.53\times10^{-5}\) within 15 DWR loops. Cells become
  increasingly stretched in the direction of the interior layer due to the
  anisotropic mesh refinement. Consequently, the maximal aspect ratio
  \(\operatorname{ar}_{\max}\) increases from 4 to 841. The effectivity index \(\Ieff\) remains bounded
  and relatively close to 1.

  Table~\ref{tab:sl-s-aniso} (structured anisotropic) also shows efficient error
  reduction. Since the structured grid is not aligned with the interior layer,
  finer meshes are required to achieve errors comparable to the unstructured
  anisotropic case. The maximal aspect ratio remains lower,
  which reflects the limitations imposed by the missing alignment of the mesh
  with the interior layer.  Overall, the anisotropic error estimators are well balanced and
  reliably track the true error while using fewer DoFs compared to the other
  strategies (cf. Tables~\ref{tab:sl-us-iso}--\ref{tab:sl-s-iso}).

  In contrast, Tables~\ref{tab:sl-us-iso}, \ref{tab:sl-s-iso}
  and~\ref{tab:sl-us-iso-uniform} (isotropic refinement on unstructured and
  structured meshes and uniform refinement) show a noticeably slower decrease in
  the error norm. The tables confirm the observations in
  Figure~\ref{fig:step-layer-error}. In all cases, the effectivity indices are
  satisfactory.

  \begin{table}[htb]\setlength{\tabcolsep}{4.5pt}
    \caption{\label{tab:sl-us-aniso}
\emph{Anisotropic} adaptive refinement on an \emph{unstructured} mesh including effectivity indices, error indicators and maximum aspect ratios for goal quantity \eqref{eq:42:JL2L2}, $\varepsilon=10^{-6}$, $\delta_0=0.1$ for Example~\ref{sec:8.1:step-layer}.
}
    \begin{center}\scriptsize
      \begin{tabular}{crrrrrrrrrrrr}
      \toprule
      $\ell$ & $N_{\text{tot}}$ &
      $N_{\text{space}}$ & $N_{\text{time}}$ &
                                                \multicolumn{2}{c}{$\|e\|_{\mathcal{Q}}$}
        & $\eta_{h,x}$  & $\eta_{h,y}$  & $\eta_h$ & $\eta_\tau$ &
                                                                   $\eta_{\tau h}$ & $\Ieff^{\textup{a}}$ & $\operatorname{ar}_{\max}$\\ \midrule
      1 & \num{1060}  & \num{53} & \num{20} & 5.333e-02 & - &                  6.143e-03 & 9.954e-03&  1.609e-02 & 2.584e-03 & 1.868e-02 & 0.35 &  4\\
      2 & \num{2220}  & \num{111} & \num{20} & 3.886e-02 & 0.46 &              4.275e-03 & 6.451e-03&  1.072e-02 & 3.798e-03 & 1.452e-02 & 0.37 &  7\\
      3 & \num{7072}  & \num{221} & \num{32} & 2.683e-02 & 0.53 &             4.144e-03 & 4.728e-03&  8.871e-03 & 1.533e-03 & 1.040e-02 & 0.39 & 14\\
      4 & \num{15488} & \num{484} & \num{32} & 1.894e-02 & 0.50 &            2.901e-03 & 2.829e-03&  5.730e-03 & 2.189e-03 & 7.919e-03 & 0.42 & 15\\
      5 & \num{44166} & \num{866} & \num{51} & 1.092e-02 & 0.79 &            1.848e-03 & 1.745e-03&  3.593e-03 & 1.078e-03 & 4.672e-03 & 0.43 & 29\\
      6 & \num{128709} & \num{1589} & \num{81} & 4.726e-03 & 1.21 &         1.147e-03 & 1.169e-03&  2.317e-03 & 9.179e-04 & 3.234e-03 & 0.68 & 30\\
      7 & \num{340689} & \num{2641} & \num{129} & 1.812e-03 & 1.38 &        2.341e-04 & 2.417e-04&  4.758e-04 & 1.116e-03 & 1.591e-03 & 0.88 & 53\\
      8 & \num{962020} & \num{4670} & \num{206} & 8.577e-04 & 1.08 &       -1.003e-04 &-1.241e-04& -2.244e-04 & 9.001e-04 & 6.757e-04 & 0.79 &105\\
      9 & \num{1536430} & \num{4670} & \num{329} & 5.790e-04 & 0.57 &       6.604e-06 &-2.690e-05& -2.030e-05 & 5.010e-04 & 4.807e-04 & 0.83 &105\\
      10 &\num{2456420} & \num{4670} & \num{526} & 4.564e-04 & 0.34 &     -1.882e-05 &-6.018e-05& -7.900e-05 & 3.099e-04 & 2.309e-04 & 0.51 &105\\
      11 &\num{3927470} & \num{4670} & \num{841} & 3.600e-04 & 0.34 &      5.735e-05 & 5.731e-06&  6.308e-05 & 1.321e-04 & 1.952e-04 & 0.54 &105\\
      12 &\num{10168200} & \num{7560} & \num{1345} & 1.559e-04 & 1.21 &   -1.777e-05 &-3.096e-05& -4.873e-05 & 1.174e-04 & 6.872e-05 & 0.44 &211\\
      13 &\num{28176136} & \num{13093} & \num{2152} & 9.794e-05 & 0.67 & -1.176e-04 &-1.232e-04& -2.408e-04 & 9.128e-05 & -1.495e-04 & 1.53&380\\
      14 &\num{76159160} & \num{22120} & \num{3443} & 6.125e-05 & 0.68 & -6.369e-05 &-7.054e-05& -1.342e-04 & 4.856e-05 & -8.567e-05 & 1.40&421\\
      15 &\num{201631356} & \num{36607} & \num{5508} & 3.534e-05 & 0.79 &-3.968e-05 &-4.136e-05& -8.105e-05 & 3.086e-05 & -5.018e-05 & 1.42&841\\
      \bottomrule
      \end{tabular}
    \end{center}
  \end{table}
  \begin{table}[htb]\setlength{\tabcolsep}{4.5pt}
    \caption{\label{tab:sl-s-aniso}\emph{Anisotropic} adaptive refinement on a \emph{structured} mesh including effectivity indices, error indicators and maximum aspect ratios for goal quantity \eqref{eq:42:JL2L2}, $\varepsilon=10^{-6}$, $\delta_0=0.1$ for Example~\ref{sec:8.1:step-layer}.}
    \begin{center}\scriptsize
      \begin{tabular}{crrrrrrrrrrrr}
        \toprule
        $\ell$ & $N_{\text{tot}}$ & $N_{\text{space}}$ & $N_{\text{time}}$ &
        \multicolumn{2}{c}{$\|e\|_{\mathcal{Q}}$}
        & $\eta_{h,x}$  & $\eta_{h,y}$  & $\eta_h$ & $\eta_\tau$ &
                                                                   $\eta_{\tau h}$ & $\Ieff^{\textup{a}}$ & $\operatorname{ar}_{\max}$\\ \midrule
        1 & \num{1620} & \num{81} & \num{20} &                 5.065e-02 &    - & 9.973e-03 &  1.651e-03& 1.162e-02 & 2.463e-03 & 1.408e-02 &  0.28 &  1 \\
        5 & \num{65856} & \num{2058} & \num{32} &            1.544e-02 & 0.46 & 5.445e-03 &  1.674e-03& 7.120e-03 & 2.607e-03 & 9.727e-03 &  0.63 &  4 \\
        6 & \num{224859} & \num{4409} & \num{51} &           9.755e-03 & 0.66 & 4.846e-03 &  1.723e-03& 6.569e-03 & 1.138e-03 & 7.708e-03 &  0.79 &  8 \\
        7 & \num{465630} & \num{9130} & \num{51} &           5.972e-03 & 0.71 & 3.254e-03 &  1.408e-03& 4.662e-03 & 1.885e-03 & 6.548e-03 &  1.10 &  8 \\
        8 & \num{1590111} & \num{19631} & \num{81} &        2.917e-03 & 1.03 & 1.539e-03 &  8.821e-04& 2.421e-03 & 1.423e-03 & 3.845e-03 &  1.32 &  8 \\
        9 & \num{5316993} & \num{41217} & \num{129} &       1.465e-03 & 0.99 & 2.420e-04 &  1.627e-04& 4.047e-04 & 1.389e-03 & 1.794e-03 &  1.23 & 16 \\
        10 &\num{15538580} & \num{75430} & \num{206} &      8.646e-04 & 0.76 & 9.463e-05 &  5.731e-05& 1.519e-04 & 8.977e-04 & 1.049e-03 &  1.21 & 16 \\
        11 &\num{24816470} & \num{75430} & \num{329} &      5.850e-04 & 0.56 & 1.774e-04 &  1.186e-04& 2.960e-04 & 4.899e-04 & 7.860e-04 &  1.34 & 16 \\
        12 &\num{73151872} & \num{139072} & \num{526} &   3.571e-04 & 0.71 &-1.352e-05 & -2.024e-05& -3.376e-05 & 3.872e-04 & 3.534e-04 &  0.99 & 16 \\
        13 &\num{116959552} & \num{139072} & \num{841} &  2.265e-04 & 0.66 & 4.424e-06 & -7.214e-06& -2.790e-06 & 2.097e-04 & 2.069e-04 &  0.91 & 16 \\
        14 &\num{187051840} & \num{139072} & \num{1345} &  1.649e-04 & 0.46 & 4.296e-05 &  2.593e-05& 6.889e-05 & 1.076e-04 & 1.765e-04 &  1.07 & 16 \\
        \bottomrule
      \end{tabular}
    \end{center}
  \end{table}
  \begin{table}[htb]\setlength{\tabcolsep}{5pt}
    \caption{\label{tab:sl-us-iso}
    \emph{Isotropic} adaptive refinement on an \emph{unstructured} mesh including effectivity indices and error indicators for goal quantity \eqref{eq:42:JL2L2}, $\varepsilon=10^{-6}$, $\delta_0=0.1$ for Example~\ref{sec:8.1:step-layer}.}
    \begin{center}\scriptsize
      \begin{tabular}{crrrrrrrrrrrr}
      \toprule
      $\ell$ & $N_{\text{tot}}$ & $N_{\text{space}}$ & $N_{\text{time}}$ &
\multicolumn{2}{c}{$\|e\|_{(0,T)\times\Omega}$} & $\eta_h$ & $\eta_\tau$ & $\eta_{\tau h}$ & $\Ieff$\\
\midrule
1 & \num{1060} & \num{53} & \num{20} & 5.3333e-02 & - & 1.9008e-02 & 2.5843e-03 & 2.1592e-02 &  0.405\\
20 & \num{348413} & \num{3819} & \num{129} & 4.3748e-03 & 0.43 & 2.1016e-03 & 3.8237e-04 & 2.4840e-03 & 0.568\\
21 & \num{405133} & \num{4427} & \num{129} & 3.9994e-03 & 0.13 & 1.6850e-03 & 4.2606e-04 & 2.1111e-03 & 0.528\\
22 & \num{465675} & \num{5195} & \num{129} & 3.6932e-03 & 0.11 & 1.2406e-03 & 4.5936e-04 & 1.7000e-03 & 0.460\\
26 & \num{1813955} & \num{9791} & \num{329} & 1.4684e-03 & 0.40 & 8.1342e-04 & 1.7237e-04 & 9.8579e-04 & 0.671\\
27 & \num{2105527} & \num{11315} & \num{329} & 1.2275e-03 & 0.26 & 7.6551e-04 & 2.2961e-04 & 9.9513e-04 & 0.811\\
28 & \num{4025176} & \num{13299} & \num{526} & 1.1618e-03 & 0.08 & 5.8804e-04 & 1.0736e-04 & 6.9539e-04 & 0.599\\
33 & \num{11337819} & \num{29493} & \num{841} & 4.8914e-04 & 0.33 & 2.2268e-04 & 9.1019e-05 & 3.1370e-04 & 0.641\\
34 & \num{20654385} & \num{33035} & \num{1345} & 3.6996e-04 & 0.40 & 2.4201e-04 & 4.3546e-05 & 2.8555e-04 & 0.772\\
35 & \num{23761175} & \num{43143} & \num{1345} & 3.3944e-04 & 0.12 & 2.1303e-04 & 4.8844e-05 & 2.6188e-04 & 0.772\\
      \bottomrule
      \end{tabular}
    \end{center}
  \end{table}
  \begin{table}[htb]\setlength{\tabcolsep}{5pt}
    \caption{\label{tab:sl-us-iso-uniform}
    \emph{Global space-time} mesh refinement on an \emph{unstructured} mesh including effectivity indices and error indicators for goal quantity \eqref{eq:42:JL2L2}, $\varepsilon=10^{-6}$, $\delta_0=0.1$ for Example~\ref{sec:8.1:step-layer}.}
    \begin{center}\scriptsize
      \begin{tabular}{crrrrrrrrrrrr}
        \toprule
        $\ell$ & $N_{\text{tot}}$ & $N_{\text{space}}$ & $N_{\text{time}}$ & \multicolumn{2}{c}{$\|e\|_{(0,T)\times\Omega}$} & $\eta_h$ & $\eta_\tau$ & $\eta_{\tau h}$ & $\Ieff$\\ \midrule
1 & \num{1060} & \num{53} & \num{20} & 5.3333e-02 & - & 1.9008e-02 & 2.5843e-03 & 2.1592e-02 & 0.405\\
2 & \num{3700} & \num{185} & \num{20} & 3.8446e-02 & 0.47 & 1.3415e-02 & 3.9320e-03 & 1.7347e-02 & 0.451\\
3 & \num{27560} & \num{689} & \num{40} & 2.6465e-02 & 0.54 & 9.4822e-03 & 1.6755e-03 & 1.1158e-02 & 0.422\\
4 & \num{106280} & \num{2657} & \num{40} & 1.8854e-02 & 0.49 & 6.1267e-03 & 2.3485e-03 & 8.4752e-03 & 0.450\\
5 & \num{834640} & \num{10433} & \num{80} & 1.0852e-02 & 0.80 & 3.9685e-03 & 1.0812e-03 & 5.0497e-03 & 0.465\\
6 & \num{3307600} & \num{41345} & \num{80} & 5.2902e-03 & 1.04 & 2.1262e-03 & 2.1816e-03 & 4.3078e-03 & 0.814\\
7 & \num{26337440} & \num{164609} & \num{160} & 1.9767e-03 & 1.42 & 6.1838e-04 & 1.4855e-03 & 2.1039e-03 & 1.064\\
8 & \num{210207040} & \num{656897} & \num{320} & 8.2674e-04 & 1.26 & 8.6590e-05 & 8.9817e-04 & 9.8476e-04 & 1.191\\
        \bottomrule
      \end{tabular}
    \end{center}
  \end{table}
  \begin{table}[htb]\setlength{\tabcolsep}{5pt}
    \caption{\label{tab:sl-s-iso}
    \emph{Isotropic} adaptive refinement on an \emph{structured} mesh including effectivity indices and error indicators for goal quantity \eqref{eq:42:JL2L2}, $\varepsilon=10^{-6}$, $\delta_0=0.1$ for Example~\ref{sec:8.1:step-layer}.}
    \begin{center}\scriptsize
      \begin{tabular}{crrrrrrrrrrrrr}
      \toprule
      $\ell$ & $N_{\text{tot}}$ & $N_{\text{space}}$ & $N_{\text{time}}$ &
                                                                            \multicolumn{2}{c}{$\|e\|_{(0,T)\times\Omega}$} & $\eta_h$ & $\eta_\tau$ & $\eta_{\tau h}$ & $\Ieff$\\
\midrule
1 & \num{1620} & \num{81} & \num{20} & 5.0651e-02 & - & 1.1831e-02 & 2.4635e-03 & 1.4294e-02 &  0.282\\
24 & \num{809029} & \num{13067} & \num{81} & 4.7925e-03 & 0.20 & 2.9181e-03 & 8.1319e-04 & 3.7313e-03 & 0.779\\
25 & \num{957281} & \num{15391} & \num{81} & 4.2042e-03 & 0.19 & 2.3924e-03 & 9.4366e-04 & 3.3361e-03 & 0.794\\
26 & \num{1716843} & \num{16939} & \num{129} & 3.3112e-03 & 0.34 & 2.0590e-03 & 5.7717e-04 & 2.6362e-03 &0.796\\
27 & \num{1971045} & \num{21699} & \num{129} & 2.7787e-03 & 0.25 & 1.7044e-03 & 6.9268e-04 & 2.3971e-03 &0.863\\
28 & \num{3681232} & \num{27849} & \num{206} & 2.2060e-03 & 0.33 & 1.6903e-03 & 3.2951e-04 & 2.0198e-03 &0.916\\
29 & \num{4481050} & \num{31829} & \num{206} & 1.9394e-03 & 0.19 & 1.4083e-03 & 3.7049e-04 & 1.7788e-03 &0.917\\
30 & \num{5302584} & \num{34857} & \num{206} & 1.6684e-03 & 0.22 & 1.0610e-03 & 4.3793e-04 & 1.4989e-03 &0.898\\
34 & \num{39173901} & \num{76007} & \num{841} & 6.2703e-04 & 0.41 & 3.8103e-04 & 6.7582e-05 & 4.4862e-04 & 0.715\\
35 & \num{47231743} & \num{90559} & \num{841} & 5.0324e-04 & 0.32 & 2.7943e-04 & 8.5870e-05 & 3.6530e-04 & 0.726\\
36 & \num{91245257} & \num{104407} & \num{1345} & 3.7876e-04 & 0.41 & 2.1721e-04 & 4.0112e-05 & 2.5732e-04 & 0.679\\
      \bottomrule
      \end{tabular}
    \end{center}
  \end{table}
  \FloatBarrier%
  \subsection{Stationary Hemker Problem}
  \label{sec:8.2:Hemker}
\begin{figure}[htb]
    \centering
    \begin{minipage}{0.48\textwidth}\centering
      \begin{tikzpicture}[scale=0.5,thick]
        \draw[black] (3,-3) -- ++(0,6);%
        \draw[black] (-3, 3) -- ++ (6,-6);
        \draw[black] (-3, -3) -- ++ (6,6);
        \draw[black] (-3, 0) -- ++ (11,0);
        \draw[black] (0, -3) -- ++ (0,6);
        \draw[blue] (-3, -3) -- ++ (11,0);
        \draw[blue] (8, -3) -- ++ (0,6);
        \draw[blue] (-3, 3) -- ++ (11,0);
        \draw[red,-latex] (2.5, 1.5) --node[above,red]{$\mathbf{b}$} ++ (2,0);
        \draw[red,-latex] (2.5, -1.5) -- ++ (2,0);
        \draw[green,fill=white] (0,0) circle (1);%
        \draw[green] (-3,-3) -- node[left]{\color{green}$\Gamma_D$}++(0,6);%
      \end{tikzpicture}
    \end{minipage}
    \begin{minipage}{0.48\textwidth}\centering
      \includegraphics[width=0.9\linewidth]{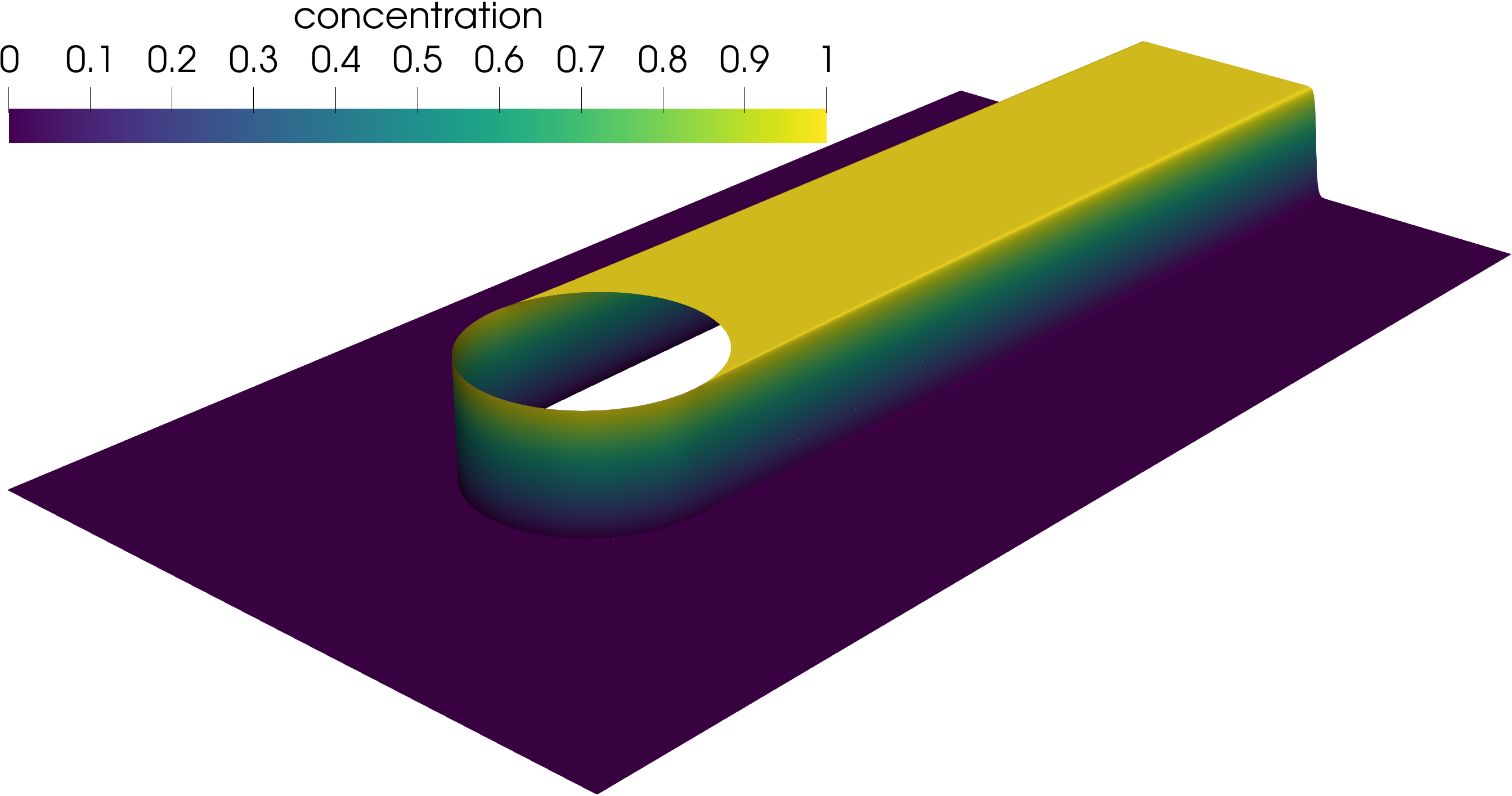}
    \end{minipage}
    \caption{\label{fig:hemker-geo}Geometry and coarse mesh of the domain for
      the Hemker problem (left) and the best adaptive solution obtained in this
      work (right). In the sketch of the geometry, green coloring corresponds to
      Dirichlet BCs, while blue indicates homogeneous Neumann BCs. On the
      circle, inhomogeneous Dirichlet BCSs are prescribed. The left boundary is
      associated with homogeneous Dirichlet BCs.}
    \includegraphics[width=0.9\linewidth]{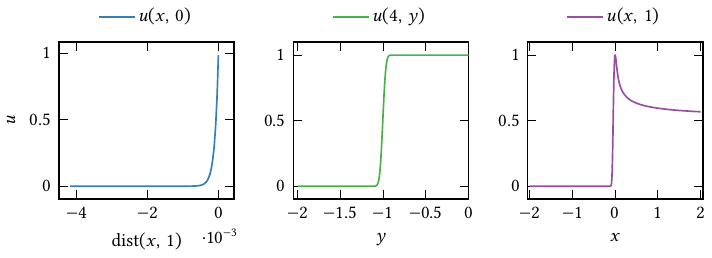}
    \caption{\label{fig:hemker-cut}Cut lines of the solution to the Hemker problem. On the left we
      plot the boundary layer in front of the obstacle. In the middle we plot a
      cut through one of the interior layers. On the right we plot the solution
      tangential to the obstacle and within the interior layer.}
  \end{figure}

\begin{figure}
    \centering
    \includegraphics[width=0.49\textwidth]{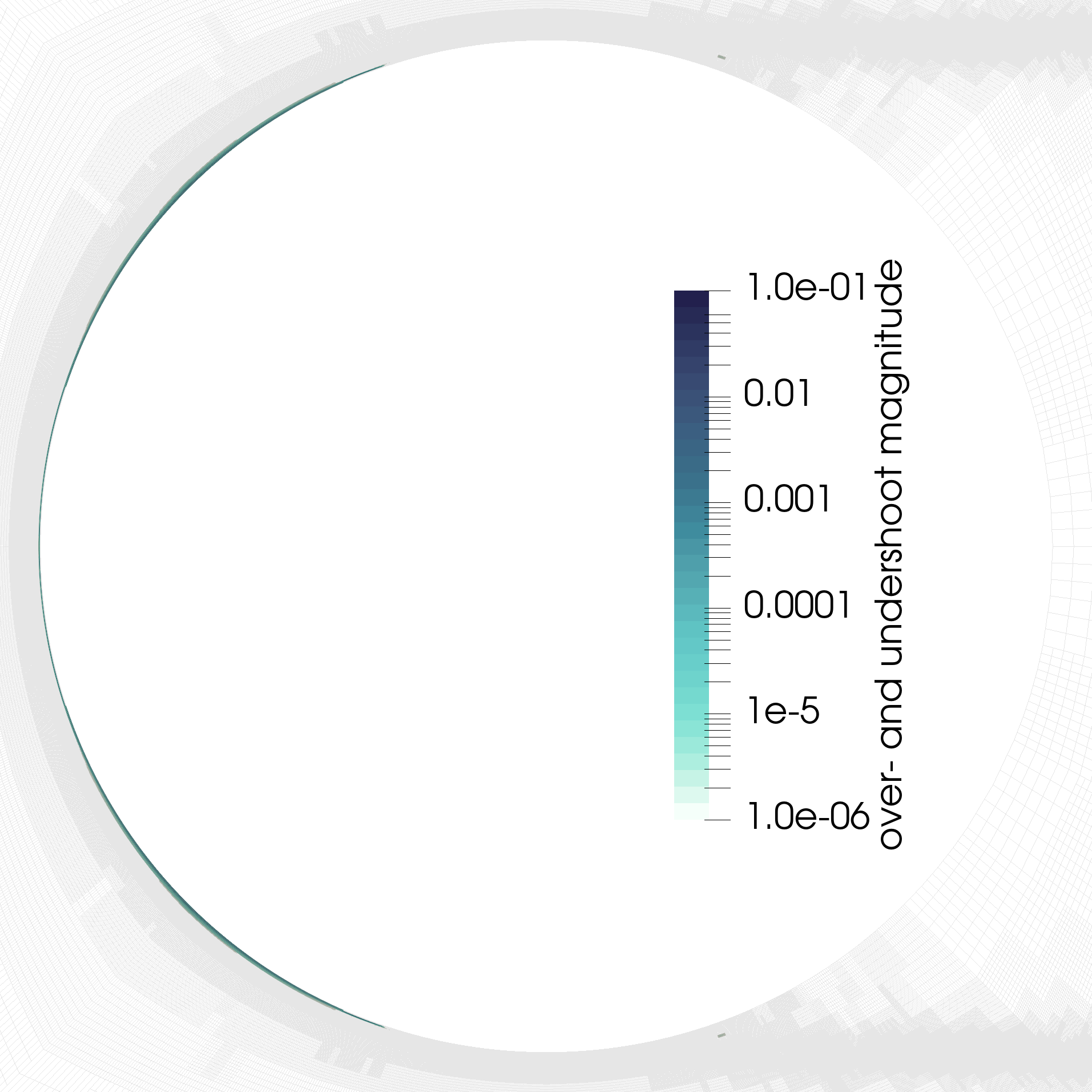}
    \includegraphics[width=0.49\textwidth]{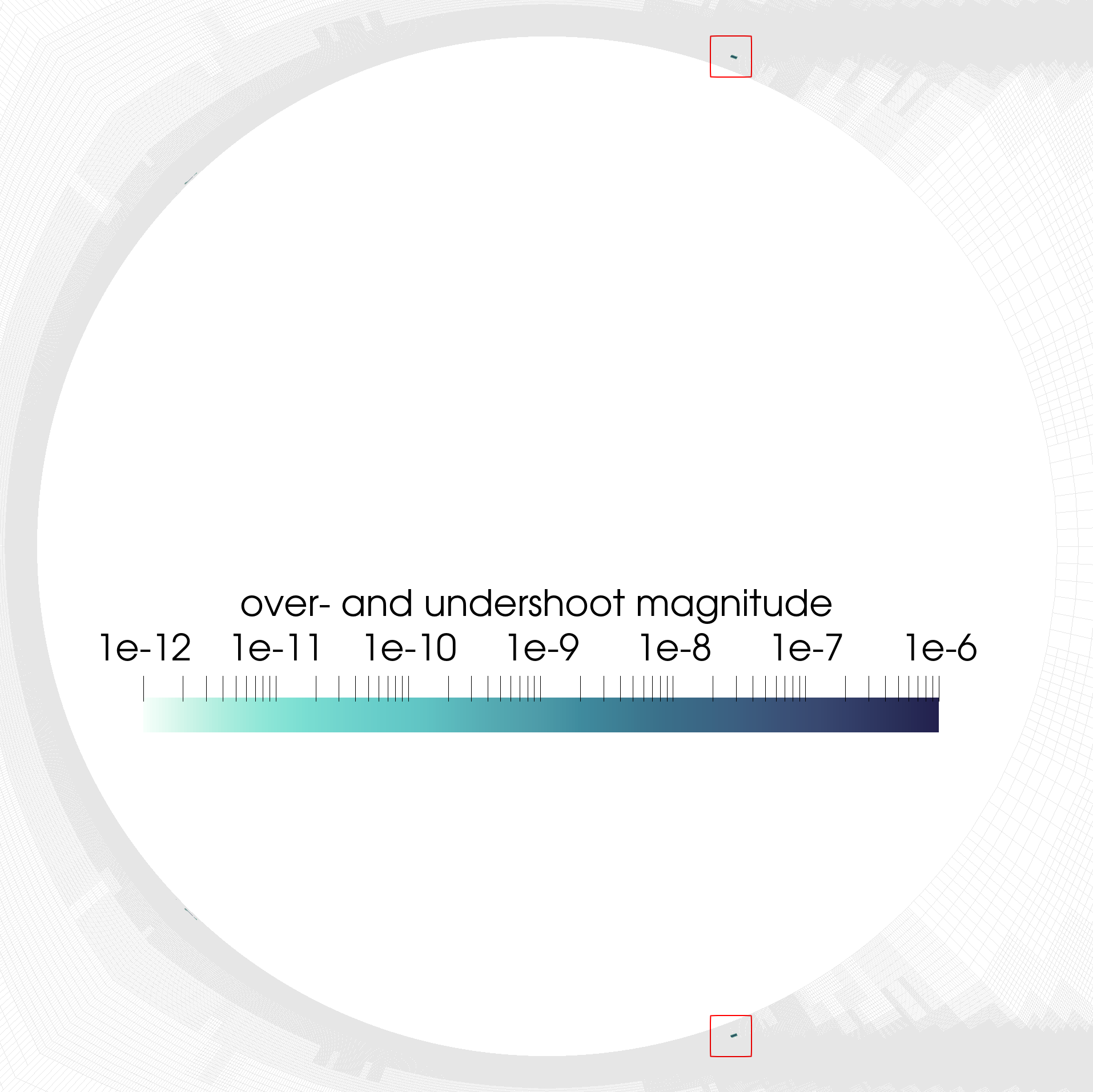}
    \caption{\label{fig:ou}The magnitude of the over- and undershoots in the DWR loop 10
      (left) and 11 (right). The maximal magnitude decreases from
      $\num{7.978e-2}$ to $\num{2.164e-6}$. Only regions where the over- and
      undershoot is larger than $\num{1e-6}$ are plotted.}
  \end{figure}

  We study the classical Hemker problem introduced
  in~\cite{hemker_singularly_1996}. It consists of the convection diffusion
  system~\eqref{eq:1:CDRoriginal} with \(\alpha = 0\), \(f = 0\) and without the
  time derivative. Here, we set $\varepsilon = \num{1.e-4}$. The concentration is entirely advected through the domain
  without any inherent degradation of the species concentration. Consequently,
  any decrease in concentration arises solely from numerical artifacts, aside
  from a negligible reduction caused by the diffusion term. The computational
  domain is defined as
  \( \Omega = \left( (-3,\,8) \times (-3,\,3) \right) \setminus \{ (x,\,y) \vert
  x^2 + y^2 \leq 1 \}\), illustrated with the coarse mesh in
  Figure~\ref{fig:hemker-geo}. The boundaries colored green in
  Figure~\ref{fig:hemker-geo} correspond to Dirichlet boundary conditions.
  Specifically, the boundary at \(x = -3\) is subjected to homogeneous Dirichlet
  conditions $u=0$. The boundaries colored blue indicate homogeneous Neumann
  boundary conditions $\partial_{\mathbf{n}} u=0$. On the circular boundary,
  inhomogeneous Dirichlet boundary conditions are imposed with \(u = 1\). The
  refinement fraction is set to
  \(\theta_{\text{space}}^{\text{ref}} = \frac{1}{3}\), and the coarsening
  fraction is \(\theta_{\text{space}}^{\text{co}} = 0\).
  Here, we use the following goal functional given by
\begin{equation}
\label{eq:45:JMean}
J(u)= \displaystyle\int_\Omega u\,\mathrm{d}\boldsymbol x\,.
\end{equation}
  The presence of sharp boundary and interior layers is further substantiated in
  Figure~\ref{fig:hemker-cut}, where the left plot shows the boundary layer in
  front of the obstacle, the middle plot provides a cut through an interior
  layer, and the right plot illustrates the solution tangential to the obstacle
  and within the interior layer.
  All cut lines confirm the correct resolution of the layers within the solution
  profile with significantly reduced oscillations.
  Table~\ref{tab:error} quantifies the evolution
  of the spatial error indicators, over successive DWR loops, showing a
  consistent reduction in error magnitude and an increase in the maximum aspect
  ratio from approximately 3.8 to 14.3. Figure~\ref{fig:ou} compares the over-
  and undershoot magnitudes between DWR loops 10 and 11, demonstrating a
  dramatic decrease in the maximal undershoot from \(\num{7.978e-2}\) to
  \(\num{2.164e-6}\). We emphasize that the undershoot of \(\num{2.164e-6}\) is
  achieved at only \(\num{3061518}\) DoFs. In the $12$th DWR loop
  maximal over- and undershoot magnitude is further reduced to $\num{3.338e-8}$.
  The width of the interior layer,
  \begin{equation}\label{eq:ylayer}
    y_{\text{layer}}=y_1-y_0\,,
  \end{equation}
  is defined to be the length of the interval $[y_0,\,y_1]$ in which the
  solution falls from $u(4, y_0) = 0.9$ to $u(4, y_1) = 0.1$. John et
  al.~\cite{augustin_assessment_2011} provide a reference value of
  $y_{\text{layer}}=0.0723$, which we achieve exactly from DWR loop 11 on (cf.\ Table~\ref{tab:error}).
  Overall, these results confirm that the goal-oriented adaptive strategy
  effectively resolves the critical features of the solution, such as boundary
  and interior layers, while simultaneously minimizing numerical artifacts.
  \begin{table}[htb]
    \caption{\label{tab:error}
    Anisotropic adaptive refinement including directional error indicators,
    maximum aspect ratios and the layer width over the DWR loops for Example~\ref{sec:8.2:Hemker}, $\varepsilon=10^{-4}$, $\delta_0=0.1$.}
    \begin{center}\footnotesize
      \begin{tabular}{crrrrrr}
      \toprule
        $\ell$ & $N_{\text{tot}}$ & $\eta_{h,x}$ & $\eta_{h,y}$ & $\eta_{h}$&
                                                                              $\operatorname{ar}_{\max}$ &$y_{\text{layer}}$\\
\midrule
        1&   \num{   414}&   -3.47e-01 & -2.09e-01 & -5.59e-01  & 3.81 &\num{0.7547}\\
        2&   \num{  1412}&   -1.56e-01 &  7.83e-02 & -7.77e-02 &  3.97 &\num{0.5997}\\
        3&   \num{  3480}&   -8.70e-02 &  1.38e-02 & -7.31e-02 &  5.95 &\num{0.5747}\\
        4&   \num{  8780}&   -2.11e-01 & -1.51e-01 & -3.61e-01 &  5.97 &\num{0.3612}\\
        5&   \num{ 22534}&    4.58e-02 &  7.20e-02 &  1.18e-01 &  6.81 &\num{0.2017}\\
        6&   \num{ 53386}&   -1.01e-02 & -2.65e-03 & -1.27e-02 &  7.52 &\num{0.1154}\\
        7&   \num{121062}&    1.85e-03 &  2.53e-03 &  4.39e-03 & 14.25 &\num{0.0883}\\
        8&   \num{272362}&   -1.25e-04 & -3.26e-04 & -4.52e-04 & 14.25 &\num{0.0720}\\
        9&   \num{621254}&   -2.26e-04 & -3.24e-04 & -5.50e-04 & 14.25 &\num{0.0727}\\
        10&  \num{1380502}&  -1.77e-04 & -2.30e-04 & -4.06e-04 & 14.25 &\num{0.0725}\\
        11&  \num{3061518}&  -1.21e-04 & -1.42e-04 & -2.63e-04 & 14.27 &\num{0.0723}\\
        12&  \num{7054500}&  -5.98e-05 & -5.95e-05 & -1.19e-04 & 14.27 &\num{0.0723}\\
      \bottomrule
      \end{tabular}
    \end{center}
  \end{table}

  \FloatBarrier%
  \subsection{Nonstationary Hemker Problem with Quadratic Obstacle}
  \label{sec:8.3:Hemker}
  In a third example, we study a modified Hemker problem by a quadratic obstacle in a
  nonstationary setting. Here, in constrast to Sec.~\ref{sec:8.2:Hemker}, we set
  $\varepsilon = \num{1.e-6}$. The computational domain is defined as
  \( \Omega = \left( (-3,\,8) \times (-3,\,3) \right) \setminus \{ (x,\,y) \vert
  \max(x,\,y) \leq 1 \}\). On the quadratic obstacle, inhomogeneous Dirichlet
  boundary conditions are imposed with \(u = 1\). The refinement fractions are
  set to \(\theta_{\text{space}}^{\text{ref}} = \frac{1}{6}\),
  \(\theta_{\text{time}}^{\text{ref}} = \frac{1}{10}\), and the coarsening
  fraction is \(\theta_{\text{space}}^{\text{co}} = 0\). The aim is to control
  the error within a control point $\boldsymbol x_{e}=(4,\,1)$, which is located within the upper interior layer. Thus, we use the
  goal functional
  \begin{equation}\label{eq:point-goal}
    J(u)=u(\boldsymbol x_{e}).
  \end{equation}
  We regularize $J$ by a
  regularized Dirac delta function
  $\delta_{r,\boldsymbol c}(\boldsymbol
  x)=\alpha\operatorname{e}^{(1-1/(1-r^2/s^2))}$, where
  $r=\lVert \boldsymbol x-\boldsymbol c\rVert$, $s>0$ is the cutoff radius and
  $\alpha$ the scaling factor, such that $\delta_{r,\boldsymbol c}$ integrates
  to $1$. This test case is solved using the $\Q{2}$ elements for the
  primal and $\Q{4}$ elements for the adjoint problem.

\begin{figure}[htbp]
    \centering
    \includegraphics[width=0.7\textwidth]{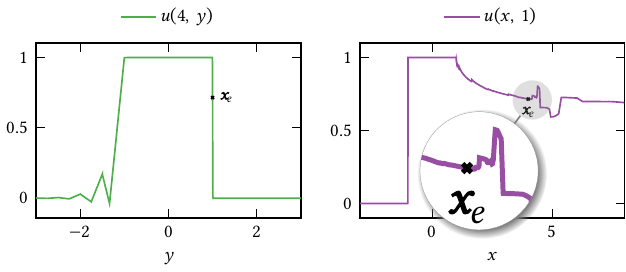}
    \caption{Cut lines of the solution to the nonstationary Hemker problem at
      final time. On the left, we
      plot a cut through the interior layers. On the right, we plot a
      cut tangential to the obstacle and within the upper interior layer.
      The point of interest is marked in both plots, around which the smoothness
      of the solution can be observed. In the right plot, we
      observe that downstream of $\boldsymbol{x}_{e}$, the solution
      deteriorates. This degradation occurs because the downstream solution has no influence on $u(\boldsymbol{x}_{e})$, and thus on the goal~\eqref{eq:point-goal}.}\label{fig:quad-hemker-cuts}
  \end{figure}

\begin{figure}[htb]
    \centering
    \begin{tikzpicture}
      \node[anchor=south west,inner sep=0, outer sep=0] (img) at (0,0)
      {\includegraphics[width=0.85\textwidth]{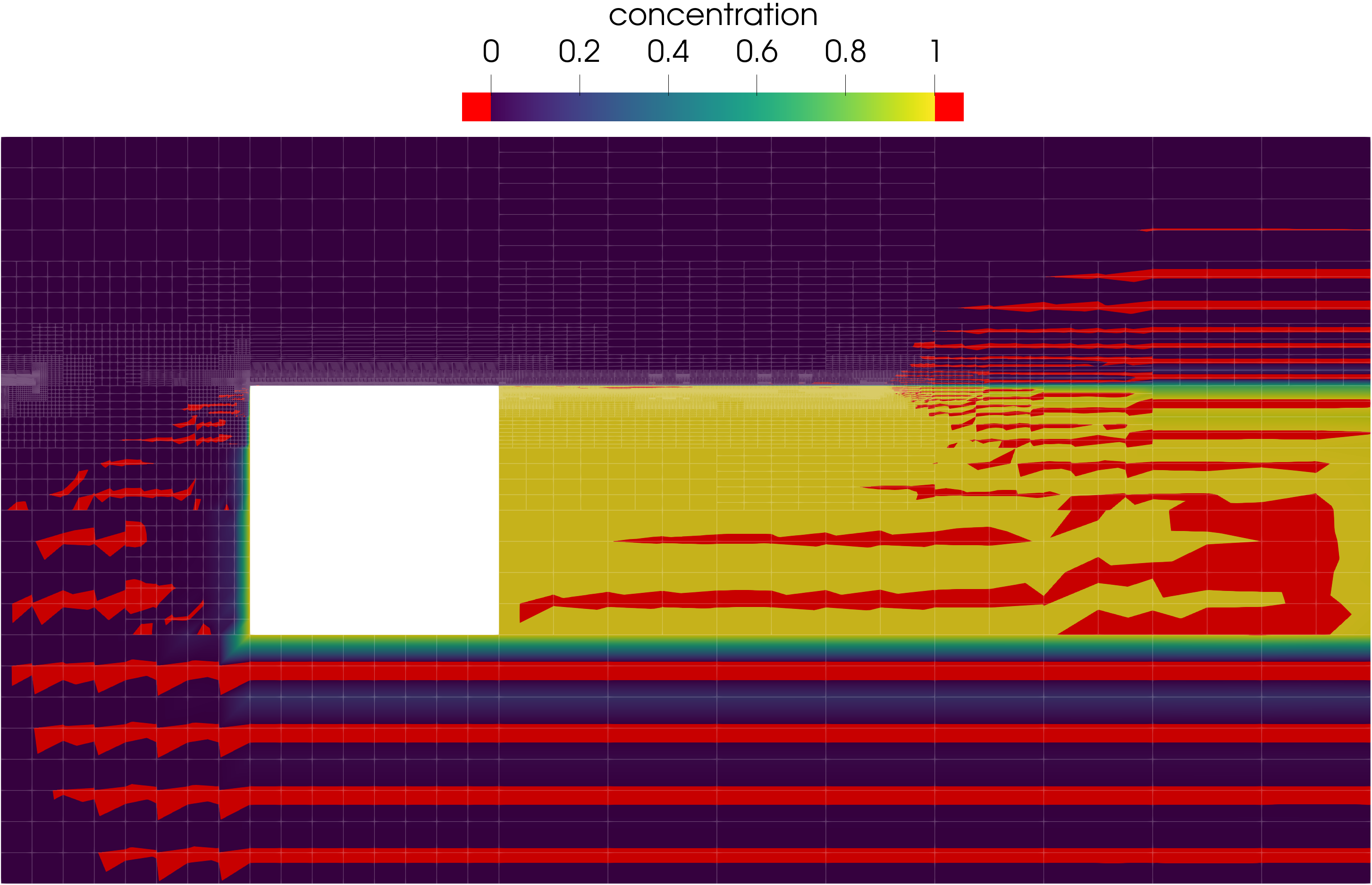}};
      \node[anchor=south,inner sep=0, outer sep=0, font=\small] at ($(img.south west) + 0.85*(img.north west)$) {$(-3,\,3)$};
      \node[anchor=south,inner sep=0, outer sep=0, font=\small] at ($(img.south east) + 0.85*(img.north west)$) {$(8,\,3)$};
      \node[anchor=north,yshift=-2pt,inner sep=0, outer sep=0, font=\small] at (img.south west) {$(-3,\,-3)$};
      \node[anchor=north,yshift=-2pt,inner sep=0, outer sep=0, font=\small] at (img.south east) {$(8,\,-3)$};
      \coordinate (pos) at ($(img.south west) + 0.6363*(img.south east) + 0.565*(img.north west)$);
      \fill[orange] (pos) circle (3pt);
      \node[above,text=orange] at (pos) {$\boldsymbol x_{e}$};
    \end{tikzpicture}
    \begin{tikzpicture}
      \node[anchor=south west,inner sep=0, outer sep=0] (img) at (0,0) {\includegraphics[width=0.85\textwidth]{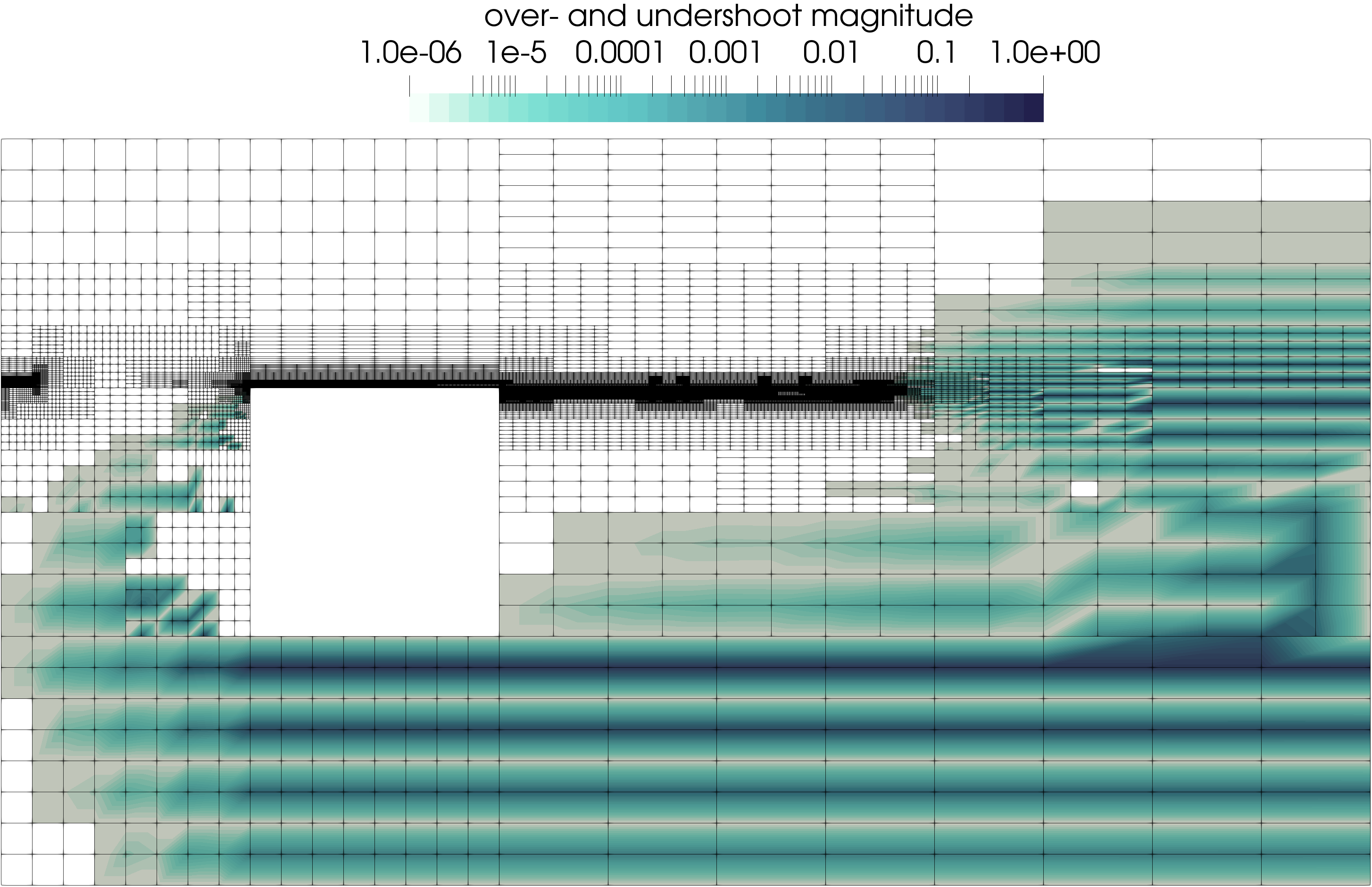}};
      \node[anchor=south,inner sep=0, outer sep=0, font=\small] at ($(img.south west) + 0.85*(img.north west)$) {$(-3,\,3)$};
      \node[anchor=south,inner sep=0, outer sep=0, font=\small] at ($(img.south east) + 0.85*(img.north west)$) {$(8,\,3)$};
      \node[anchor=north,yshift=-2pt,inner sep=0, outer sep=0, font=\small] at (img.south west) {$(-3,\,-3)$};
      \node[anchor=north,yshift=-2pt,inner sep=0, outer sep=0, font=\small] at (img.south east) {$(8,\,-3)$};
      \coordinate (pos) at ($(img.south west) + 0.6363*(img.south east) + 0.565*(img.north west)$);
      \fill[orange] (pos) circle (3pt);
      \node[left,yshift=-2pt,text=orange] at (pos) {$\boldsymbol x_{e}$};
    \end{tikzpicture}
    \caption{The solution and over- and undershoot magnitude in the final DWR
      loop for Example~\ref{sec:8.3:Hemker}. One can clearly see the adaptive mesh refinement at
      $\boldsymbol{x}_{e}$ and upstream of the point due to the goal~\eqref{eq:point-goal}.}\label{fig:quad-hemker-solution}
  \vspace*{3ex}
  \includegraphics[width=0.85\textwidth]{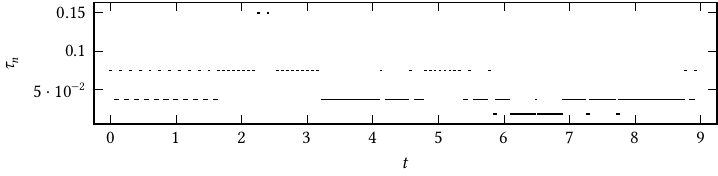}
  \vspace*{-2ex}
    \caption{The distribution of timestep sizes over the
      interval $I$.}\label{fig:quad-hemker-tsteps}
  \end{figure}

  The cut-line plots in Figure~\ref{fig:quad-hemker-cuts} provide insight into
  the effectiveness of the anisotropic adaptive mesh refinement. The plots confirm the
  smoothness around \(\boldsymbol{x}_e\) and display a degradation downstream,
  which is consistent with the goal-oriented error control at the control point.
  Figure~\ref{fig:quad-hemker-solution} illustrates the solution and its over-
  and undershoot magnitudes. We observe the adaptive mesh refinement in
  the vicinity of \(\boldsymbol{x}_e\) and regions upstream of
  \(\boldsymbol{x}_e\). In Figure~\ref{fig:quad-hemker-tsteps}, the distribution
  of timestep sizes shows that the smallest timesteps occur for
  \( t \in [6,7] \), coinciding with the passage of the solution front through
  the control point \(\boldsymbol{x}_e\). Moreover,
  Table~\ref{tab:quad-hemker-estimators} quantifies the evolution of the spatial
  and temporal error estimators over the DWR loops. The convergence of the error
  estimators can be clearly observed. The maximum aspect ratio of the spatial
  mesh increases from 3.5 to 256, underscoring the effectiveness of anisotropic
  refinement.
  \begin{table}[htb]\setlength{\tabcolsep}{5pt}
    \caption{
      Anisotropic adaptive refinement including directional error indicators for goal quantity~\eqref{eq:point-goal} and maximum aspect ratios for Example~\ref{sec:8.3:Hemker}, $\varepsilon=10^{-6}$, $\delta_0=0.1$.}\label{tab:quad-hemker-estimators}
    \centering
    \scriptsize
    \begin{tabular}{crrrrrrrrrr}
    \toprule
    $\ell$ & $N_{\text{tot}}$ &
    $N_{\text{space}}$ & $N_{\text{time}}$ &
    $\eta_{h,x}$  & $\eta_{h,y}$  & $\eta_h$ & $\eta_\tau$ &
    $\eta_{\tau h}$ & $\operatorname{ar}_{\max}$\\
    \midrule
       1& \num{   576}&  60&  \num{   34560}&   1.8381e-24& -2.7453e-21& -2.7434e-21& -4.7170e-23& -2.7906e-21&     3.5 \\
       2& \num{   914}&  66&  \num{   60324}&  -3.0163e-03& -6.4458e-01& -6.4760e-01& -7.6284e-03& -6.5522e-01&     3.5 \\
       3& \num{  1594}&  72&  \num{  114768}&  -2.1019e-03& -7.6246e-01& -7.6456e-01& -6.7205e-03& -7.7128e-01&     4.0 \\
       4& \num{  2996}&  79&  \num{  236684}&  -1.7760e-03& -2.6215e-01& -2.6392e-01& -1.8290e-03& -2.6575e-01&     8.0 \\
       5& \num{  4934}&  86&  \num{  424324}&  -1.6904e-03& -1.1922e-01& -1.2091e-01& -7.0710e-04& -1.2162e-01&    16.0 \\
       6& \num{  8672}&  94&  \num{  815168}&  -1.4765e-03& -4.4182e-02& -4.5658e-02& -2.1388e-04& -4.5872e-02&    16.0 \\
       7& \num{ 14642}& 103&  \num{ 1508126}&  -1.0071e-03& -9.2494e-03& -1.0256e-02& -4.1194e-05& -1.0298e-02&    16.0 \\
       8& \num{ 25914}& 113&  \num{ 2928282}&  -5.2393e-04& -6.2303e-04& -1.1470e-03& -1.0574e-05& -1.1575e-03&    32.0 \\
       9& \num{ 44614}& 124&  \num{ 5532136}&  -2.9816e-04& -1.2836e-04& -4.2652e-04& -1.3834e-05& -4.4036e-04&    64.0 \\
      10& \num{ 73420}& 136&  \num{ 9985120}&  -1.7129e-04& -6.9224e-05& -2.4051e-04& -1.5069e-05& -2.5558e-04&   128.0 \\
      11& \num{112842}& 149&  \num{16813458}&  -4.8920e-05& -3.9307e-06& -5.2850e-05& -1.5029e-05& -6.7879e-05&   256.0 \\
      \bottomrule
    \end{tabular}
    \end{table}


  \FloatBarrier%
  \section{Conclusion \& Outlook}
  In this work, we proposed an anisotropic goal-oriented error estimator based
  on the DWR method for time-dependent CDR equations.
  The error estimator decouples directional error contributions in space and
  time while incurring only negligible additional cost compared to the standard
  DWR approach. Moreover, these error indicators quantify anisotropy of the
  solution with respect to the goal. SUPG stabilization was applied to suppress
  spurious oscillations at high Péclet numbers.

  The proposed anisotropic adaptive algorithm produces high aspect ratio
  elements. Consequently, the numerical solutions capture sharp layers without
  producing spurious oscillations. Numerical experiments confirm that
  goal-oriented anisotropic refinement in combination with SUPG stabilization
  outperforms isotropic and global strategies in both accuracy and computational
  efficiency. The underlying algorithm yields efficient and robust results for
  multiple challenging benchmarks including both interior and boundary
  layers. The implementation can address complex geometries by supporting
  unstructured meshes with curved elements. Therefore, our proposed method
  effectively controls the error while maintaining feasible problem sizes in the
  convection-dominated regime.

  Currently the implemenation only supports tensor-products of a fixed
  spatial mesh with the temporal mesh. Future work includes the extension to
  dynamic meshes~\cite{BruchhaeuserSV08} and quasilinear parabolic problems~\cite{TOULOPOULOS2025127}. Furthermore, an extension to $hp$
  refinement and discontinuous Galerkin discretizations in space is intended.
  The methods presented here can be applied to other problems in a
  straightforward manner.

  \section*{Acknowledgments}
  All authors acknowledge the funding of DAAD-project 57729992,
  "Goal-oriented AnIsotropic Space-Time Mesh Adaption (AIMASIM)" in the funding program "Programm des projektbezogenen Personenaustauschs Griechenland ab 2024".
  Bernhard Endtmayer and Thomas Wick additionally acknowledge the support by the Cluster of Excellence PhoenixD (EXC 2122, Project ID 390833453).
  Bernhard Endtmayer was funded by an Humboldt Postdoctoral Fellowship at the beginning of the work.
  Bernhard Endtmayer thanks the research group of Markus Bause for financing the
  research visit at the Helmut Schmidt University Hamburg in February 16--28,
  2025.
  Computational resources (HPC cluster HSUper) have been provided by the project
  hpc.bw, funded by dtec.bw - Digitalization and Technology Research Center of the
  Bundeswehr. dtec.bw is funded by the European Union - NextGenerationEU.


  \appendix
\section*{Appendix}
The explicit definitions of the adjoint bilinear forms in \eqref{eq:3:9:A_tau_prime_u_phi_z_eq_J_prime_u_phi_A_S_prime_u_phi_z_eq_J_prime_u_phi} are given by
\begin{align*}
  A^{\prime}(u)(\varphi,z)
  & =  \int_{I}\big\{
(\varphi,- \partial_{t} z)
+ a^{\prime}(u)(\varphi,z)
\big\}\mathrm{d}t + (\varphi(T),z(T))\,,
\\
  A_{\tau}^{\prime}(u_{\tau})(\varphi_{\tau},z_{\tau})
  & = \sum_{n=1}^N\int_{I_n}\big\{
(\varphi_{\tau},-\partial_{t} z_{\tau})
+ a^{\prime}(u_{\tau})(\varphi_{\tau},z_{\tau})
\big\}\mathrm{d}t
 - \sum_{n=1}^{N-1}(\varphi_{\tau,n}^-,[z_{\tau}]_{n})
+(\varphi_{\tau,N}^-, z_{\tau,N}^-)
    \,,
  \\
  A_{S}^{\prime}(u_{\tau h})(\varphi_{\tau h},z_{\tau h})
  & =  \sum_{n=1}^N\int_{I_n}\big\{
(\varphi_{\tau h},- \partial_{t} z_{\tau h})
+ a_{h}^{\prime}(u_{\tau h})(\varphi_{\tau h},z_{\tau h})
\big\}\mathrm{d}t
\\
& \;\;\;\,
+ S^{\prime}(u_{\tau h})(\varphi_{\tau h},z_{\tau h})
- \sum_{n=1}^{N-1}(\varphi_{\tau h,n}^-,[z_{\tau h}]_{n})
+ (\varphi_{\tau h,N}^-, z_{\tau h,N}^-)
    \,,
\end{align*}
where the derivative $a^{\prime}(u)(\varphi,z)$ of the
bilinear form $a(u)(z)$ admits
the explicit form
\begin{displaymath}
a^{\prime}(u)(\varphi,z)=
(\varepsilon\nabla \varphi,\nabla z)
+ (\boldsymbol{b} \cdot \nabla \varphi,z)
+ (\alpha \varphi,z)\,.
\end{displaymath}
We note that for the representation of
$A^{\prime}(\cdot, \cdot)(\cdot, \cdot)$ integration by parts in time is applied,
which is allowed for weak
solutions $z\in \mathcal{V}$; cf., e.g.,~\cite[Lemma 8.9]{BruchhaeuserR17}.
%
\bibliographystyle{abbrv}
\bibliography{lit}

\end{document}

%% file: header.tex
\usepackage{amsmath,amsthm,amsfonts,amssymb,amscd}
\usepackage{graphicx}
\usepackage{mathtools}
\usepackage{hyperref}
\usepackage{pgf,tikz}
\usepackage{pgfplots}
\usepackage{tikz-3dplot}
\usepackage{mathrsfs}
\usepackage{float}
\usepackage{csvsimple}
\usepackage{enumerate}
\usepackage[inline]{enumitem}
\usepackage{booktabs}
\usepackage{pgfplots}
\pgfplotsset{compat=newest}
\usepackage{graphicx,bm,xcolor}
\usepackage{algorithm}
\usepackage{pdfpages}
\usepackage[noend]{algorithmic}
\renewcommand{\algorithmiccomment}[1]{\bgroup\hfill$\blacktriangleright$~#1\egroup}
\usepackage{stmaryrd}
\usepackage[numbers]{natbib}
\usepackage{hyperref}
\usetikzlibrary{arrows}
\usepackage{caption}
\usepackage[T1]{fontenc}
\usepackage{t1enc}
\usepackage{verbatim} 
\usepackage{xcolor}
\usepackage{subcaption}
\usepackage{upgreek}

\usepackage{url}

\input{Tikz/Function_for_drawings.tex}
\newif\ifMAKEPICS
\MAKEPICSfalse

\ifMAKEPICS
\usepackage[cleanup,subfolder]{gnuplottex}
\usepackage{xparse}

\ExplSyntaxOn
\DeclareExpandableDocumentCommand{\convertlen}{ O{cm} m }
{
	\dim_to_decimal_in_unit:nn { #2 } { 1 #1 } cm
}
\ExplSyntaxOff
\fi


%% file: Tikz/Function_for_drawings.tex
\newcommand{\functionforplot}[2]{(#1^2 - 2*#1 + #2^2/2 + #1*#2/4 -#2 + 1)}

%% file: Tikz/FigureQ1.tex
\definecolor{zzttqq}{rgb}{0.2,0.3,0.4}
\definecolor{xdxdff}{rgb}{0.1,0.1,0.5}
\definecolor{ududff}{rgb}{0.8,0.8,0}
\definecolor{uuuuuu}{rgb}{0.4,0.4,0.4}
\begin{tikzpicture}[line cap=round,line join=round,>=triangle 45,x=2cm,y=2cm]
\clip(-0.4,-0.4) rectangle (2.4,2.4);
\fill[line width=2pt,color=zzttqq,fill=blue,fill opacity=0.2] (0,0) -- (2,0) -- (2,2) -- (0,2) -- cycle;
\fill[line width=2pt,color=zzttqq,fill=zzttqq,fill opacity=0.0] (0,1) -- (1,1) -- (1,2) -- (0,2) -- cycle;
\fill[line width=2pt,color=zzttqq,fill=zzttqq,fill opacity=0.0] (1,1) -- (2,1) -- (2,2) -- (1,2) -- cycle;
\fill[line width=2pt,color=zzttqq,fill=zzttqq,fill opacity=0.0] (0,0) -- (1,0) -- (1,1) -- (0,1) -- cycle;
\fill[line width=2pt,color=zzttqq,fill=zzttqq,fill opacity=0.0] (1,0) -- (2,0) -- (2,1) -- (1,1) -- cycle;
\draw [line width=2pt] (0,2)-- (2,2);
\draw [line width=2pt] (2,2)-- (2,0);
\draw [line width=2pt] (2,0)-- (0,0);
\draw [line width=2pt] (0,0)-- (0,2);
\draw [line width=2pt] (0,1)-- (2,1);
\draw [line width=2pt] (1,2)-- (1,0);
\draw [line width=2pt,color=uuuuuu] (0,0)-- (2,0);
\draw [line width=2pt,color=uuuuuu] (2,0)-- (2,2);
\draw [line width=2pt,color=uuuuuu] (2,2)-- (0,2);
\draw [line width=2pt,color=uuuuuu] (0,2)-- (0,0);
\draw [line width=2pt,color=uuuuuu] (0,1)-- (1,1);
\draw [line width=2pt,color=uuuuuu] (1,1)-- (1,2);
\draw [line width=2pt,color=uuuuuu] (1,2)-- (0,2);
\draw [line width=2pt,color=uuuuuu] (0,2)-- (0,1);
\draw [line width=2pt,color=uuuuuu] (1,1)-- (2,1);
\draw [line width=2pt,color=uuuuuu] (2,1)-- (2,2);
\draw [line width=2pt,color=uuuuuu] (2,2)-- (1,2);
\draw [line width=2pt,color=uuuuuu] (1,2)-- (1,1);
\draw [line width=2pt,color=uuuuuu] (0,0)-- (1,0);
\draw [line width=2pt,color=uuuuuu] (1,0)-- (1,1);
\draw [line width=2pt,color=uuuuuu] (1,1)-- (0,1);
\draw [line width=2pt,color=uuuuuu] (0,1)-- (0,0);
\draw [line width=2pt,color=uuuuuu] (1,0)-- (2,0);
\draw [line width=2pt,color=uuuuuu] (2,0)-- (2,1);
\draw [line width=2pt,color=uuuuuu] (2,1)-- (1,1);
\draw [line width=2pt,color=uuuuuu] (1,1)-- (1,0);
\begin{scriptsize}
\draw [fill=xdxdff] (0,0) circle (2.5pt);
\draw[color=uuuuuu] (0.0,-0.15) node {$N_1$};
\draw [fill=xdxdff] (2,2) circle (2.5pt);
\draw[color=uuuuuu] (2.0,2.15) node {$N_9$};
\draw [fill=xdxdff] (0,2) circle (2.5pt);
\draw[color=uuuuuu] (0.0,2.15) node {$N_7$};
\draw [fill=xdxdff] (2,0) circle (2.5pt);
\draw[color=uuuuuu] (2.0,-0.15) node {$N_3$};
\draw [fill=xdxdff] (1,1) circle (2.5pt);
\draw[color=uuuuuu] (1.15,1.15) node {$N_5$};
\draw [fill=xdxdff] (0,1) circle (2.5pt);
\draw[color=uuuuuu] (-0.15,1.0) node {$N_4$};
\draw [fill=xdxdff] (2,1) circle (2.5pt);
\draw[color=uuuuuu] (2.14,1.0) node {$N_6$};
\draw [fill=xdxdff] (1,2) circle (2.5pt);
\draw[color=uuuuuu] (1.0,2.15) node {$N_8$};
\draw[color=uuuuuu] (1.0,-0.15) node {$N_2$};
\draw [fill=xdxdff] (1,0) circle (2.5pt);
\draw [fill=xdxdff] (2,2) circle (2.5pt);
\draw [fill=xdxdff] (0,2) circle (2.5pt);
\draw [fill=xdxdff] (1,2) circle (2.5pt);
\draw[color=uuuuuu] (1.5,1.5) node {$K_{h,4}$};
\draw[color=uuuuuu] (0.5,1.5) node {$K_{h,3}$};
\draw[color=uuuuuu] (1.5,0.5) node {$K_{h,2}$};
\draw[color=uuuuuu] (0.5,0.5) node {$K_{h,1}$};
\draw [fill=xdxdff] (0,2) circle (2.5pt);
\draw [fill=xdxdff] (2,2) circle (2.5pt);
\draw [fill=xdxdff] (1,2) circle (2.5pt);
\draw [fill=xdxdff] (1,1) circle (2.5pt);
\draw [fill=xdxdff] (0,1) circle (2.5pt);
\draw [fill=xdxdff] (2,1) circle (2.5pt);
\draw [fill=xdxdff] (1,1) circle (2.5pt);
\end{scriptsize}
\end{tikzpicture}

%% file: Tikz/FigureQ2.tex
\definecolor{zzttqq}{rgb}{0.2,0.3,0.4}
\definecolor{xdxdff}{rgb}{0.1,0.1,0.5}
\definecolor{ududff}{rgb}{0.8,0.8,0}
\definecolor{uuuuuu}{rgb}{0.4,0.4,0.4}
\begin{tikzpicture}[line cap=round,line join=round,>=triangle 45,x=2cm,y=2cm]
\clip(-0.4,-0.4) rectangle (2.4,2.4);
\fill[line width=2pt,color=green,fill=orange,fill opacity=0.2] (0,0) -- (2,0) -- (2,2) -- (0,2) -- cycle;
\draw [line width=2pt] (0,2)-- (2,2);
\draw [line width=2pt] (2,2)-- (2,0);
\draw [line width=2pt] (2,0)-- (0,0);
\draw [line width=2pt] (0,0)-- (0,2);
\begin{scriptsize}
\draw [fill=xdxdff] (0,0) circle (2.5pt);
\draw [fill=xdxdff] (2,2) circle (2.5pt);
\draw [fill=xdxdff] (0,2) circle (2.5pt);
\draw [fill=xdxdff] (2,0) circle (2.5pt);
\draw [fill=xdxdff] (1,1) circle (2.5pt);
\draw [fill=xdxdff] (0,1) circle (2.5pt);
\draw [fill=xdxdff] (2,1) circle (2.5pt);
\draw [fill=xdxdff] (1,2) circle (2.5pt);
\draw [fill=xdxdff] (1,0) circle (2.5pt);
\draw [fill=xdxdff] (2,2) circle (2.5pt);
\draw [fill=xdxdff] (0,2) circle (2.5pt);
\draw [fill=xdxdff] (1,2) circle (2.5pt);
\draw[color=uuuuuu] (0.75,0.75) node {$Q_2$};
\draw [fill=xdxdff] (0,2) circle (2.5pt);
\draw [fill=xdxdff] (2,2) circle (2.5pt);
\draw [fill=xdxdff] (1,2) circle (2.5pt);
\draw [fill=xdxdff] (1,1) circle (2.5pt);
\draw [fill=xdxdff] (0,1) circle (2.5pt);
\draw [fill=xdxdff] (2,1) circle (2.5pt);
\draw [fill=xdxdff] (1,1) circle (2.5pt);

\draw[color=uuuuuu] (0.0,-0.15) node {$N_1$};
\draw [fill=xdxdff] (2,2) circle (2.5pt);
\draw[color=uuuuuu] (2.0,2.15) node {$N_9$};
\draw [fill=xdxdff] (0,2) circle (2.5pt);
\draw[color=uuuuuu] (0.0,2.15) node {$N_7$};
\draw [fill=xdxdff] (2,0) circle (2.5pt);
\draw[color=uuuuuu] (2.0,-0.15) node {$N_3$};
\draw [fill=xdxdff] (1,1) circle (2.5pt);
\draw[color=uuuuuu] (1.15,1.15) node {$N_5$};
\draw [fill=xdxdff] (0,1) circle (2.5pt);
\draw[color=uuuuuu] (-0.15,1.0) node {$N_4$};
\draw [fill=xdxdff] (2,1) circle (2.5pt);
\draw[color=uuuuuu] (2.14,1.0) node {$N_6$};
\draw [fill=xdxdff] (1,2) circle (2.5pt);
\draw[color=uuuuuu] (1.0,2.15) node {$N_8$};
\draw[color=uuuuuu] (1.0,-0.15) node {$N_2$};

\end{scriptsize}
\end{tikzpicture}

%% file: Tikz/PatchQ11.tex
\begin{tikzpicture}
\begin{axis}[
xlabel={$x$},                  
ylabel={$y$},               
hide z axis,
grid=both,                    
samples=3,                   
domain=0:2,                   
y domain=0:2,                 
xtick={0,1,2},
xticklabels={},
yticklabels={},
xlabel={$x_1$},
ylabel={$x_2$},
colormap/jet,
mesh/ordering=x varies,
mesh/cols=13
]
\addplot3[
surf,                         
shader=interp,                  
opacity=0.3                   
]
{\functionforplot{x}{y}};  

\addplot3[
mesh,
draw=black,
samples=3, 
]
{\functionforplot{x}{y}};  
\end{axis}
\end{tikzpicture}

%% file: Tikz/PatchQ22HighInter.tex
\begin{tikzpicture}
\begin{axis}[
xlabel={$x$},                  
ylabel={$y$},               
hide z axis,
grid=both,                    
samples=30,                   
domain=0:2,                   
y domain=0:2,                 
xtick={0,1,2},
xticklabels={},
yticklabels={},
xlabel={$x_1$},
ylabel={$x_2$},
colormap/jet,
mesh/ordering=x varies,
mesh/cols=13
]
\addplot3[
surf,                         
shader=interp,                  
opacity=0.3                   
]
{\functionforplot{x}{y}};  

\addplot3 [domain=0:2,samples y=1] (1,x,\functionforplot{1}{x});
\addplot3 [domain=0:2,samples y=1] (0,x,\functionforplot{0}{x});
\addplot3 [domain=0:2,samples y=1] (2,x,\functionforplot{2}{x});
\addplot3 [domain=0:2,samples y=1] (x,1,\functionforplot{x}{1});
\addplot3 [domain=0:2,samples y=1] (x,0,\functionforplot{x}{0});
\addplot3 [domain=0:2,samples y=1] (x,2,\functionforplot{x}{2});
\end{axis}
\end{tikzpicture}

%% file: Tikz/FigureQ21.tex
\definecolor{zzttqq}{rgb}{0.2,0.3,0.4}
\definecolor{xdxdff}{rgb}{0.1,0.1,0.5}
\definecolor{ududff}{rgb}{0.8,0.8,0}
\definecolor{uuuuuu}{rgb}{0.4,0.4,0.4}
\begin{tikzpicture}[line cap=round,line join=round,>=triangle 45,x=2cm,y=2cm]
	\clip(-0.4,-0.4) rectangle (2.4,2.4);
	\fill[line width=2pt,color=zzttqq,fill=orange,fill opacity=0.15] (0,0) -- (2,0) -- (2,2) -- (0,2) -- cycle;
	\fill[line width=2pt,color=zzttqq,fill=zzttqq,fill opacity=0.0] (0,1) -- (1,1) -- (1,2) -- (0,2) -- cycle;
	\fill[line width=2pt,color=zzttqq,fill=zzttqq,fill opacity=0.0] (1,1) -- (2,1) -- (2,2) -- (1,2) -- cycle;
	\fill[line width=2pt,color=zzttqq,fill=zzttqq,fill opacity=0.0] (0,0) -- (1,0) -- (1,1) -- (0,1) -- cycle;
	\fill[line width=2pt,color=zzttqq,fill=zzttqq,fill opacity=0.0] (1,0) -- (2,0) -- (2,1) -- (1,1) -- cycle;
	\draw [line width=2pt] (0,2)-- (2,2);
	\draw [line width=2pt] (2,2)-- (2,0);
	\draw [line width=2pt] (2,0)-- (0,0);
	\draw [line width=2pt] (0,0)-- (0,2);
	\draw [line width=2pt] (0,1)-- (2,1);
	\begin{scriptsize}
		\draw [fill=xdxdff] (0,0) circle (2.5pt);
		\draw[color=uuuuuu] (0.0,-0.15) node {$N_1$};
		\draw [fill=xdxdff] (2,2) circle (2.5pt);
		\draw[color=uuuuuu] (2.0,2.15) node {$N_9$};
		\draw [fill=xdxdff] (0,2) circle (2.5pt);
		\draw[color=uuuuuu] (0.0,2.15) node {$N_7$};
		\draw [fill=xdxdff] (2,0) circle (2.5pt);
		\draw[color=uuuuuu] (2.0,-0.15) node {$N_3$};
		\draw [fill=xdxdff] (1,1) circle (2.5pt);
		\draw[color=uuuuuu] (1.15,1.15) node {$N_5$};
		\draw [fill=xdxdff] (0,1) circle (2.5pt);
		\draw[color=uuuuuu] (-0.15,1.0) node {$N_4$};
		\draw [fill=xdxdff] (2,1) circle (2.5pt);
		\draw[color=uuuuuu] (2.14,1.0) node {$N_6$};
		\draw [fill=xdxdff] (1,2) circle (2.5pt);
		\draw[color=uuuuuu] (1.0,2.15) node {$N_8$};
		\draw[color=uuuuuu] (1.0,-0.15) node {$N_2$};
		\draw [fill=xdxdff] (1,0) circle (2.5pt);
		\draw [fill=xdxdff] (2,2) circle (2.5pt);
		\draw [fill=xdxdff] (0,2) circle (2.5pt);
		\draw [fill=xdxdff] (1,2) circle (2.5pt);
		\draw[color=uuuuuu] (1,1.5) node {$K_{h,3,4}$};
		\draw[color=uuuuuu] (1,0.5) node {$K_{h,1,2}$};
		\draw [fill=xdxdff] (0,2) circle (2.5pt);
		\draw [fill=xdxdff] (2,2) circle (2.5pt);
		\draw [fill=xdxdff] (1,2) circle (2.5pt);
		\draw [fill=xdxdff] (1,1) circle (2.5pt);
		\draw [fill=xdxdff] (0,1) circle (2.5pt);
		\draw [fill=xdxdff] (2,1) circle (2.5pt);
		\draw [fill=xdxdff] (1,1) circle (2.5pt);
	\end{scriptsize}
\end{tikzpicture}

%% file: Tikz/FigureQ12.tex
\definecolor{zzttqq}{rgb}{0.2,0.3,0.4}
\definecolor{xdxdff}{rgb}{0.1,0.1,0.5}
\definecolor{ududff}{rgb}{0.8,0.8,0}
\definecolor{uuuuuu}{rgb}{0.4,0.4,0.4}
\begin{tikzpicture}[line cap=round,line join=round,>=triangle 45,x=2cm,y=2cm]
\clip(-0.4,-0.4) rectangle (2.4,2.4);
\fill[line width=2pt,color=zzttqq,fill=green,fill opacity=0.15] (0,0) -- (2,0) -- (2,2) -- (0,2) -- cycle;
\fill[line width=2pt,color=zzttqq,fill=zzttqq,fill opacity=0.0] (0,1) -- (1,1) -- (1,2) -- (0,2) -- cycle;
\fill[line width=2pt,color=zzttqq,fill=zzttqq,fill opacity=0.0] (1,1) -- (2,1) -- (2,2) -- (1,2) -- cycle;
\fill[line width=2pt,color=zzttqq,fill=zzttqq,fill opacity=0.0] (0,0) -- (1,0) -- (1,1) -- (0,1) -- cycle;
\fill[line width=2pt,color=zzttqq,fill=zzttqq,fill opacity=0.0] (1,0) -- (2,0) -- (2,1) -- (1,1) -- cycle;
\draw [line width=2pt] (0,2)-- (2,2);
\draw [line width=2pt] (2,2)-- (2,0);
\draw [line width=2pt] (2,0)-- (0,0);
\draw [line width=2pt] (0,0)-- (0,2);
\draw [line width=2pt] (1,2)-- (1,0);
\begin{scriptsize}
\draw [fill=xdxdff] (0,0) circle (2.5pt);
\draw[color=uuuuuu] (0.0,-0.15) node {$N_1$};
\draw [fill=xdxdff] (2,2) circle (2.5pt);
\draw[color=uuuuuu] (2.0,2.15) node {$N_9$};
\draw [fill=xdxdff] (0,2) circle (2.5pt);
\draw[color=uuuuuu] (0.0,2.15) node {$N_7$};
\draw [fill=xdxdff] (2,0) circle (2.5pt);
\draw[color=uuuuuu] (2.0,-0.15) node {$N_3$};
\draw [fill=xdxdff] (1,1) circle (2.5pt);
\draw[color=uuuuuu] (1.15,1.15) node {$N_5$};
\draw [fill=xdxdff] (0,1) circle (2.5pt);
\draw[color=uuuuuu] (-0.15,1.0) node {$N_4$};
\draw [fill=xdxdff] (2,1) circle (2.5pt);
\draw[color=uuuuuu] (2.14,1.0) node {$N_6$};
\draw [fill=xdxdff] (1,2) circle (2.5pt);
\draw[color=uuuuuu] (1.0,2.15) node {$N_8$};
\draw[color=uuuuuu] (1.0,-0.15) node {$N_2$};
\draw [fill=xdxdff] (1,0) circle (2.5pt);
\draw [fill=xdxdff] (2,2) circle (2.5pt);
\draw [fill=xdxdff] (0,2) circle (2.5pt);
\draw [fill=xdxdff] (1,2) circle (2.5pt);
\draw[color=uuuuuu] (1.5,1) node {$K_{h,2,4}$};
\draw[color=uuuuuu] (0.5,1) node {$K_{h,1,3}$};
\draw [fill=xdxdff] (0,2) circle (2.5pt);
\draw [fill=xdxdff] (2,2) circle (2.5pt);
\draw [fill=xdxdff] (1,2) circle (2.5pt);
\draw [fill=xdxdff] (1,1) circle (2.5pt);
\draw [fill=xdxdff] (0,1) circle (2.5pt);
\draw [fill=xdxdff] (2,1) circle (2.5pt);
\draw [fill=xdxdff] (1,1) circle (2.5pt);
\end{scriptsize}
\end{tikzpicture}

%% file: Tikz/PatchQ11HighInter.tex
\begin{tikzpicture}
\begin{axis}[
xlabel={$x$},                  
ylabel={$y$},               
hide z axis,
grid=both,                    
samples=3,                   
domain=0:2,                   
y domain=0:2,                 
xtick={0,1,2},
xticklabels={},
yticklabels={},
xlabel={$x_1$},
ylabel={$x_2$},
colormap/jet,
mesh/ordering=x varies,
mesh/cols=13
]
\addplot3[
surf,                         
shader=interp,                  
opacity=0.3                   
]
{\functionforplot{x}{y}};  

\addplot3[
mesh,
draw=black,
samples=3, 
]
{\functionforplot{x}{y}};  
\end{axis}
\end{tikzpicture}

%% file: Tikz/PatchQ21HighInter.tex
\begin{tikzpicture}
\begin{axis}[
xlabel={$x$},                 
ylabel={$y$},     
hide z axis,           
grid=both,                    
samples=30,                   
samples y=30,
domain=0:2,                   
y domain=0:2,                 
xtick={0,1,2},
xticklabels={},
yticklabels={},
xlabel={$x_1$},
ylabel={$x_2$},
colormap/jet,
mesh/ordering=x varies,
mesh/cols=13
]

\addplot3[
surf,                         
shader=interp,  
samples=3, 	              
samples y=30,
opacity=0.3                   
]
{\functionforplot{x}{y}};  
\addplot3 [domain=0:2,samples y=1,color=red,ultra thick] (1,x,\functionforplot{1}{x});
\addplot3 [domain=0:2,samples y=1,color=red,ultra thick] (0,x,\functionforplot{0}{x});
\addplot3 [domain=0:2,samples y=1,color=red,ultra thick] (2,x,\functionforplot{2}{x});
\addplot3 [domain=0:2,samples =3,samples y=1] (x,1,\functionforplot{x}{1});
\addplot3 [domain=0:2,samples =3,samples y=1] (x,0,\functionforplot{x}{0});
\addplot3 [domain=0:2,samples =3,samples y=1] (x,2,\functionforplot{x}{2});
\end{axis}

\end{tikzpicture}

%% file: Tikz/PatchQ12HighInter.tex
\begin{tikzpicture}
\begin{axis}[
xlabel={$x$},                 
ylabel={$y$},     
hide z axis,           
grid=both,                    
samples=30,                   
samples y=3,
domain=0:2,                   
y domain=0:2,                 
xtick={0,1,2},
xticklabels={},
yticklabels={},
xlabel={$x_1$},
ylabel={$x_2$},
colormap/jet,
mesh/ordering=x varies,
mesh/cols=13
]

\addplot3[
surf,                         
shader=interp,  	              
samples y=3,
opacity=0.3                   
]
{\functionforplot{x}{y}};  
\addplot3 [domain=0:2,samples =3,samples y=1] (1,x,\functionforplot{1}{x});
\addplot3 [domain=0:2,samples =3,samples y=1] (0,x,\functionforplot{0}{x});
\addplot3 [domain=0:2,samples =3,samples y=1] (2,x,\functionforplot{2}{x});
\addplot3 [domain=0:2,samples y=1,color=red,ultra thick] (x,1,\functionforplot{x}{1});
\addplot3 [domain=0:2,samples y=1,color=red,ultra thick] (x,0,\functionforplot{x}{0});
\addplot3 [domain=0:2,samples y=1,color=red,ultra thick] (x,2,\functionforplot{x}{2});
\end{axis}

\end{tikzpicture}